\pdfoutput=1
\documentclass[11pt,a4paper]{article}
\usepackage[margin=27mm,headheight=15pt]{geometry}
\usepackage[T1]{fontenc}
\usepackage[utf8]{inputenc}
\usepackage[english]{babel}
\usepackage{lmodern}
\usepackage{amsmath,amssymb,amsthm,mathtools,mathrsfs}
\usepackage{microtype,booktabs,array,enumitem,fancyhdr}
\usepackage{xcolor}
\usepackage[colorlinks=true,linkcolor=blue!45!black,citecolor=blue!45!black,urlcolor=blue!45!black]{hyperref}

\hypersetup{pdftitle={Rank-Two Bundles on Ruled Threefolds over Abelian Surfaces},pdfauthor={Qi Sun}}
\setlist{nosep,leftmargin=*}
\newtheorem{introthm}{Theorem}
\newtheorem{theorem}{Theorem}[section]
\newtheorem{proposition}[theorem]{Proposition}
\newtheorem{lemma}[theorem]{Lemma}
\newtheorem{corollary}[theorem]{Corollary}
\theoremstyle{definition}

\newtheorem{remark}[theorem]{Remark}
\newtheorem{example}[theorem]{Example}
\newcommand{\OO}{\mathcal O}
\newcommand{\VV}{\mathcal V}

\newcommand{\II}{\mathcal I}
\newcommand{\PP}{\mathbb P}
\newcommand{\CC}{\mathbb C}
\newcommand{\ZZ}{\mathbb Z}
\newcommand{\QQ}{\mathbb Q}
\newcommand{\BN}{\operatorname{BN}}
\newcommand{\Pic}{\operatorname{Pic}}
\newcommand{\NS}{\operatorname{NS}}
\newcommand{\Ext}{\operatorname{Ext}}
\newcommand{\Hom}{\operatorname{Hom}}
\newcommand{\rk}{\operatorname{rk}}

\newcommand{\Sym}{\operatorname{Sym}}

\newcommand{\cHom}{\mathcal{H}\!om}
\newcommand{\cEnd}{\mathcal{E}\!nd}

\newcommand{\pr}{\operatorname{pr}}
\newcommand{\id}{\operatorname{id}}

\title{Rank-Two Bundles on Ruled Threefolds over Abelian Surfaces}
\author{Qi Sun\thanks{School of Mathematical Sciences, Shanghai Key Laboratory of PMMP, East China Normal University, Math.\ Bldg, No.~500, Dongchuan Road, Shanghai, 200241, P.~R.~China, \href{mailto:52285500003@stu.ecnu.edu.cn}{52285500003@stu.ecnu.edu.cn}}}
\date{}

\begin{document}
\maketitle
\vspace{-12mm}
\begin{abstract}
Let $X=\PP_A(\VV)$ be a ruled threefold over a complex abelian surface $A$ with $\VV$ nef. We construct smooth projective connected components of moduli spaces of rank-two bundles on $X$ with fiber degree one via line-bundle extensions. These components are projective bundles over $\widehat A\times\widehat A$, and their Brill--Noether schemes are pulled back from cohomology jump schemes on $A$. On an abelian variety of arbitrary dimension, we determine the scheme structure of the degree-zero cohomology jump loci ($h^0$-jump schemes) for any nef semihomogeneous bundle, showing that they are cylinders over finite subschemes. When $-K_X$ is nef, this yields a trichotomy governed by $c_1(\VV)$: the jump schemes are finite of four types, effective Cartier divisors of four types having codimension one whereas the expected codimension is $k^2$, or either empty or all of $\widehat A$.
\end{abstract}



\section{Introduction}\label{sec:intro}

For a family of stable sheaves, the condition $h^0(F)\ge k$ defines a Brill--Noether locus, whose determinantal scheme structure records how global sections vary in families~\cite{CMR,BW}. On abelian surfaces, the weak Brill--Noether property was characterized by Coskun--Nuer--Yoshioka~\cite{CNY}. In Picard rank one, Bayer--Li~\cite{BL} studied higher strata when the first Chern class generates the N\'eron--Severi group, and Chien~\cite{Chien} treated the first locus for positive multiples of the generator and nonpositive Euler characteristic. On ruled threefolds over $\PP^2$, Costa and Mac\'ias Tarr\'io~\cite{CMT} constructed components of stable rank-two bundles and established the nonemptiness of their Brill--Noether strata. Replacing the rational base by an abelian surface, a natural question is which Brill--Noether scheme structures arise as a rank-two extension varies in its Picard parameters.

In this paper we study this question for a natural class of rank-two bundles on $\PP^1$-bundles over abelian surfaces whose determinant has degree one on the fibers. They are extensions of two line bundles, so their moduli are governed by the two Picard parameters, and a relative computation transfers their cohomology to the abelian surface, where Fourier--Mukai theory applies.

Let $\pi:X=\PP_A(\VV)\to A$ be the projective bundle of a nef rank-two vector bundle on a complex abelian surface. We use the quotient convention $\pi_*\OO_X(1)=\VV$ and write $\zeta=c_1(\OO_X(1))$ and $B=\pi^*H$ for an ample line bundle $H$ on $A$. For an integer $a>0$ and $P_1,P_2\in\widehat A=\Pic^0(A)$, set
\[
 S(P_1):=\OO_X(-1)\otimes\pi^*(H^{\otimes a}\otimes P_1),\qquad Q(P_2):=\OO_X(2)\otimes\pi^*P_2,
\]
and consider the extensions
\begin{equation}\label{eq:introext}
 0\longrightarrow S(P_1)\longrightarrow E\longrightarrow Q(P_2)\longrightarrow 0,
\end{equation}
whose middle terms $E$ are rank-two vector bundles with numerical first Chern class $c_1(E)=\zeta+aB$, so that $\det E$ has degree one on the fibers of $\pi$. Write $c=c_1(\VV)$. By relative duality, an extension class $\xi$ in~\eqref{eq:introext} is a section of $\VV^*\otimes(\det\VV)^{-1}\otimes H^{\otimes a}\otimes P_1P_2^{-1}$ on $A$. For $\xi\ne0$, let $D\subset A$ be the divisorial part of its zero scheme (allowing $D=0$); evaluation identifies
\[
 T_\xi:=\OO_X(1)\otimes\pi^*\pi_*(E(-1))\simeq\OO_X(1)\otimes\pi^*\bigl((\det\VV)^2\otimes H^{\otimes(-a)}\otimes P_1^{-1}P_2^2\otimes\OO_A(D)\bigr)
\]
with a saturated invertible subsheaf of $E$, where $E(-1)=E\otimes\OO_X(-1)$. For rational $t>0$ set $L_t=\zeta+tB$, and write $\mu_{L_t}(F)=c_1(F)\cdot L_t^2/\rk F$ for the slope of a torsion-free sheaf $F$. For a vector bundle $F$ on an abelian variety $A$, let $Z_k(F)\subset\widehat A$ denote the cohomology jump scheme defined by $h^0(A,F\otimes P)\ge k$. Our first result describes exactly when these extensions are stable and shows that they constitute a full connected component of the moduli space.

\begin{introthm}\label{thm:introA}
Suppose $aH-3c$ is ample.
\begin{enumerate}[label=(\roman*)]
\item For a nonzero extension $E$ in~\eqref{eq:introext} of class $\xi$ and any rational $t>0$, $E$ is $\mu_{L_t}$-stable (resp.\ $\mu_{L_t}$-semistable) if and only if
\[
 \max\{\mu_{L_t}(S(P_1)),\,\mu_{L_t}(T_\xi)\}<\mu_{L_t}(E)\qquad (\text{resp.\ }\le).
\]
In particular, every nonzero extension in~\eqref{eq:introext} is $\mu_{L_t}$-stable for $t$ in the common chamber
\begin{equation}\label{eq:introchamber}
 \tfrac a3 B\cdot L_t^2<\zeta\cdot L_t^2<(aB-2\pi^*c)\cdot L_t^2,
\end{equation}
which is optimal whenever $\VV$ admits a numerically trivial quotient line bundle.
\item For $a\gg0$ and $t$ satisfying~\eqref{eq:introchamber}, the nonzero extensions~\eqref{eq:introext} form a smooth projective connected component $\mathcal D_a$ of the Gieseker moduli space of dimension $h+3$, where $h=\chi(A,\VV^*\otimes(\det\VV)^{-1}\otimes H^{\otimes a})$. Its Albanese morphism $E\mapsto(P_1,P_2)$ realizes $\mathcal D_a$ as a $\PP^{h-1}$-bundle over $\widehat A\times\widehat A$, and scheme-theoretically $\BN_X^k\cap\mathcal D_a=q_2^{-1}(Z_k(\Sym^2\VV))$, where $q_2:\mathcal D_a\to\widehat A$ is the second projection.
\item On $X=A\times\PP^1$ ($\VV=\OO_A^{\oplus2}$) with $a=1$, $\mathcal D_1$ exists on the optimal chamber $2/3<t<2$ and persists as a smooth projective connected component of the Gieseker moduli space at the upper wall $L_2=\zeta+2B$. At this wall, the strictly $\mu_{L_2}$-semistable members---corresponding to proportional pairs of sections of $H\otimes P_1P_2^{-1}$---remain Gieseker-stable.
\end{enumerate}
\end{introthm}

Parts~(i), (ii) and~(iii) of Theorem~\ref{thm:introA} are proved in Theorem~\ref{thm:nefexact} with Corollary~\ref{cor:nefoptimal}, in Theorem~\ref{thm:oddcomponent}, and in Proposition~\ref{prop:giesekerwall}, respectively. Both Picard parameters enter the extension space, but the last statement of~(ii) shows that only $P_2$ affects global sections: since $S(P_1)$ has relative degree $-1$, pushing forward~\eqref{eq:introext} along $\pi$ gives $R\pi_*E\simeq\Sym^2\VV\otimes P_2$. The Brill--Noether problem on $X$ thus reduces to the jump schemes $Z_k:=Z_k(\Sym^2\VV)$ on the base. Our second result determines such schemes in the broader setting of nef semihomogeneous bundles on an abelian variety of arbitrary dimension, which include all homogeneous bundles twisted by nef line bundles. We use the normalized Poincar\'e bundle for the Fourier--Mukai transform $\Phi_A$.

\begin{introthm}\label{thm:introB}
Let $E$ be a nef semihomogeneous vector bundle on an abelian variety $A$ of arbitrary dimension. Set $K:=K(\det E)^0$, $s:=\dim K$, and let $\rho:\widehat A\to\widehat K$ be the restriction map.
\begin{enumerate}[label=(\roman*)]
\item For every $k\ge1$, there exists a finite subscheme $D_k\subset\widehat K$ such that $Z_k(E)=\rho^{-1}(D_k)$; in particular, every nonempty proper $Z_k(E)$ has codimension $s$. When $E=\mathcal W\otimes L$ with $\mathcal W$ homogeneous and $L$ a nef line bundle, $D_k$ is cut out by $\operatorname{Fitt}_{k-1}\bigl((N^\vee)^{\oplus m}\bigr)$, where $N$ is the finite-length Fourier--Mukai transform of $E|_K$, $N^\vee$ is its local dual, and $m\ge1$ is an integer determined by $L$.
\item When $\dim A=2$, set $\delta:=c_1(E)/\rk E$. Depending on $\delta$, the scheme $Z_k(E)$ is finite if $\delta\equiv0$; an effective Cartier divisor supported on finitely many translates of an elliptic curve if $\delta\not\equiv0$ and $\delta^2=0$ (with multiplicities governed by the Atiyah decomposition of $E|_K$ when $E=\mathcal W\otimes L$); and either $\widehat A$ or empty if $\delta$ is ample.
\end{enumerate}
\end{introthm}

To apply Theorem~\ref{thm:introB}, with its explicit description of $D_k$, to $E=\Sym^2\VV$, one needs $\Sym^2\VV$ to be a line-bundle twist of a homogeneous bundle. As we show in Lemma~\ref{lem:antinef}, the rank-three bundle $\Sym^2\VV\otimes(\det\VV)^{-1}$ is homogeneous if and only if $-K_X$ is nef. In this situation Theorem~\ref{thm:introB} applies with $L=\det\VV$, and the numerical class of $c_1(\VV)$ governs the geometry of the jump schemes, yielding the following classification.

\begin{introthm}\label{thm:introC}
Assume $-K_X$ is nef and put $c=c_1(\VV)$.
\begin{enumerate}[label=(\roman*)]
\item If $c\equiv0$, each $Z_k$ is a finite subscheme of $\widehat A$ of length at most six.
\item If $c\not\equiv0$ and $c^2=0$, each nonempty $Z_k$ is an effective Cartier divisor supported on translates of an elliptic curve in $\widehat A$; its actual codimension is one, whereas its determinantal expected codimension is $k^2$.
\item If $c$ is ample, $Z_k$ is $\widehat A$ or empty.
\end{enumerate}
In cases (i) and (ii) exactly four types occur, determined by $\VV$ and by its restriction to an elliptic curve $C\subset A$ respectively; they are listed in Theorem~\ref{thm:homogeneous} and Corollary~\ref{cor:trichotomy}.
\end{introthm}

Case~(i) of Theorem~\ref{thm:introC} follows from Mukai's correspondence between homogeneous bundles and finite-length sheaves on $\widehat A$~\cite{Mukai78,Mukai81}: the transform is a length-three sheaf in degree two, and the Fitting ideals of its local dual give the jump schemes; the types include the cube of a maximal ideal and a curvilinear scheme with ideal $(y,x^3)$. The main new input is Theorem~\ref{thm:introB}. Semihomogeneous bundles were introduced and classified by Mukai~\cite{Mukai78}. Their cohomological support loci have been studied as sets, and through cohomological rank functions~\cite{AP}; here we determine the scheme structure of all $h^0$-jump loci, including nonreduced multiplicities arising from non-simple homogeneous factors. For $E=\mathcal W\otimes L$ the transform lies in degree $s$ and has the same transverse completed local model at every point of a translate of $\ker\rho$, so the jump ideals are pulled back from $\widehat K$ and the jump schemes are cylinders over Fitting schemes of a finite-length module. Mukai's structure theory reduces an arbitrary nef semihomogeneous bundle to this case after an isogeny, which preserves the cylinder structure.

The hypothesis on $-K_X$ cannot be dropped: for $\VV=\OO_A\oplus f^*M$, with $f$ an elliptic fibration of $A$, the jump schemes acquire embedded points and components of codimension two (Example~\ref{rem:mixed}).

\paragraph{Organization.} The geometry of $X$, namely stability, moduli and walls, is studied only over abelian surfaces, whereas the jump schemes are treated on abelian varieties of arbitrary dimension and then specialized to $g=2$. Section~\ref{sec:conventions} fixes conventions and introduces the extensions and their jumping fibers. Section~\ref{sec:nefexact} proves Theorem~\ref{thm:introA}: the stability criterion and its optimality, the moduli component, and the upper wall on $A\times\PP^1$. Section~\ref{sec:jump} works on an abelian variety of arbitrary dimension and proves the statements of Theorem~\ref{thm:introB} on jump schemes. Section~\ref{sec:classification} applies this to $X$ and proves Theorem~\ref{thm:introC}, with examples.

\section{Preliminaries and the extension family}\label{sec:conventions}

This section fixes notation and conventions, introduces the extensions studied in this paper, and describes how their splitting type varies over $A$.

\subsection{Conventions}

We work over $\CC$. Throughout, $A$ is a complex abelian surface, $H$ an ample line bundle on $A$, $d=H^2$, and $\widehat A=\Pic^0(A)$; no assumption is made on the Picard rank. Let $\mathcal P$ be the \emph{normalized Poincar\'e line bundle} on $A\times\widehat A$~\cite{Mumford,BiL}, uniquely characterized up to isomorphism by the properties that $\mathcal P|_{A\times\{[P]\}}\simeq P$ for every $[P]\in\widehat A$ and $\mathcal P|_{\{0_A\}\times\widehat A}\simeq\OO_{\widehat A}$, where $0_A\in A$ is the origin. We use $H^{\otimes a}=\OO_A(aH)$ for tensor powers and $H^2$ for the self-intersection number. We write $\equiv$ for numerical equivalence. Let $\pi:X=\PP_A(\VV)\to A$ be the projective bundle of a rank-two vector bundle on $A$, in the quotient convention $\pi_*\OO_X(1)=\VV$. Write $\zeta=c_1(\OO_X(1))$. For rational $t>0$ set $L_t=\zeta+tB$ (where $B=\pi^*H$ is introduced below) and denote slope with respect to it by $\mu_{L_t}$. Then $K_X=-2\zeta+\pi^*c_1(\VV)$, and
\begin{equation}\label{eq:proj}
 R\pi_*\OO_X=\OO_A,\quad
 R\pi_*\OO_X(j)=\Sym^j\VV\ (j\ge0),\quad
 R\pi_*\OO_X(-1)=0.
\end{equation}
Twisting $\VV$ by a line bundle $N$ leaves $X$ unchanged, replaces $\zeta$ by $\zeta+\pi^*c_1(N)$, and leaves $-K_X$ unchanged.

On $X$, $M_{X,L}(w)$ denotes the \emph{stable locus} of the Gieseker moduli space with a fixed full numerical Chern character $w$ (for a rational ample class $L=L_t$, Gieseker stability is defined using any sufficiently divisible integral multiple of $L$). We do not fix a determinant unless stated.

For a flat family whose fiber cohomology vanishes in degrees at least two, the relative cohomology on a local universal chart is represented by a two-term complex of locally free sheaves $K^0\xrightarrow{d^0}K^1$. The Brill--Noether scheme is defined by
\begin{equation}\label{eq:BNideal}
 \BN^k=V\bigl(I_{\rk K^0-k+1}(d^0)\bigr).
\end{equation}
This is independent of the chosen presentation and descends from charts; it is the usual determinantal convention~\cite{CMR}. We use reduced closure when naming an irreducible component obtained from a parameter family.

When higher cohomology does not vanish, as in Section~\ref{sec:classification}, we define the scheme by the degree-zero cohomology jump ideal of a local perfect complex $K^0\to K^1\to K^2\to\cdots$ starting in degree zero. It is still $I_{\rk K^0-k+1}(d^0)$, now interpreted as a cohomology jump ideal~\cite[Section~2, especially Corollaries~2.3--2.5]{BW}; the two-term expected-codimension formula need not apply.

\subsection{The extension family}

We now introduce the extension family studied throughout the paper. Because the line subbundle and the quotient line bundle have fiber degrees $-1$ and $2$, each extension class is a global section of a twisted rank-two bundle on $A$. The pointwise vanishing of this section determines the splitting type of the middle term on the fibers, while its divisorial zero locus controls slope stability in Section~\ref{sec:nefexact}. We assume throughout that $\VV$ is nef, equivalently that $\zeta$ is nef. This normalization is always available, since a sufficiently ample twist of $\VV$ is globally generated and leaves $X$ unchanged. Put $B=\pi^*H$ and
\[
 d=H^2,\quad u=c_1(\VV)\cdot H,\quad
 v_0=\int_A(c_1(\VV)^2-c_2(\VV)),\quad
 W=\VV^*\otimes(\det\VV)^{-1}.
\]
Then $\zeta$ and $B$ are nef, $\zeta+B$ is integral ample, and
\[
 \zeta B^2=d,\quad \zeta^2B=u,\quad\zeta^3=v_0.
\]
Fix a positive integer $a$ sufficiently large that, for every $P\in\widehat A$,
\begin{equation}\label{eq:oddvanish}
 \begin{split}
 &H^i(A,W(aH)\otimes P)=0\quad(i>0),\qquad
 h:=\chi(A,W(aH))>0,\\
 &H^i(A,\Sym^3\VV(-aH)\otimes P)=0\quad(i=0,1),
 \end{split}
\end{equation}
and also
\begin{equation}\label{eq:oddachoice}
 aH-3c_1(\VV)\ \text{is ample},\qquad
 a^2d+5au+3v_0>0,\qquad
 a^2d>5au+2\textstyle\int_Ac_1(\VV)^2+v_0.
\end{equation}
Such an integer exists: relative Serre vanishing on $A\times\widehat A$ gives the first line of \eqref{eq:oddvanish} uniformly, and Serre duality gives the second line from the positive-twist vanishings of $(\Sym^3\VV)^*$. The Euler characteristic and the quadratic expressions in \eqref{eq:oddachoice} are eventually positive.

For $P_1,P_2\in\widehat A$, consider
\begin{equation}\label{eq:oddextension}
 0\longrightarrow S(P_1)\longrightarrow E\longrightarrow Q(P_2)\longrightarrow 0,
\end{equation}
recalling that $S(P_1)=\OO_X(-1)\otimes\pi^*(H^{\otimes a}\otimes P_1)$ and $Q(P_2)=\OO_X(2)\otimes\pi^*P_2$.

\begin{lemma}\label{lem:oddExt}
For every pair $(P_1,P_2)\in\widehat A\times\widehat A$, the extension group $\Ext_X^1(Q(P_2),S(P_1))$ parametrizing~\eqref{eq:oddextension} has dimension $h$, and the remaining $\Ext$-groups satisfy
\[
 \Ext_X^i(Q(P_2),S(P_1))=0\quad(i\ne1),\qquad
 \Ext_X^i(S(P_1),Q(P_2))=0\quad(i=0,1,3).
\]
Every middle term $E$ in~\eqref{eq:oddextension} is locally free of rank two, with numerical Chern classes
\[
 c_1(E)=\zeta+aB,\quad c_2(E)=-2\zeta^2+2a\zeta B,\quad c_3(E)=0.
\]
In particular, $\det E$ has degree one on the fibers of $\pi$, and for any non-split extension~\eqref{eq:oddextension}, the restriction of $E$ to a general fiber of $\pi$ splits as $\OO_{\PP^1}\oplus\OO_{\PP^1}(1)$.
\end{lemma}

\begin{proof}
In the quotient convention, the relative canonical bundle is $\omega_{X/A}=\OO_X(-2)\otimes\pi^*\det\VV$. Since $R^1\pi_*\omega_{X/A}\simeq\OO_A$ by relative duality, the projection formula gives $R^1\pi_*\OO_X(-2)\simeq(\det\VV)^{-1}$. The natural multiplication pairing
\[
 \pi_*\OO_X(1)\otimes R^1\pi_*\OO_X(-3)\longrightarrow R^1\pi_*\OO_X(-2)
\]
restricts on each fiber to the perfect Serre duality pairing $H^0(\PP^1,\OO(1))\otimes H^1(\PP^1,\OO(-3))\to H^1(\PP^1,\OO(-2))$. Since $\pi_*\OO_X(1)=\VV$, this pairing induces an isomorphism
\[
 R^1\pi_*\OO_X(-3)\simeq\mathcal{H}om\bigl(\VV,(\det\VV)^{-1}\bigr)=\VV^*\otimes(\det\VV)^{-1}=W,
\]
while $\pi_*\OO_X(-3)=0$. Consequently, the Leray spectral sequence degenerates at $E_2$ and yields
\[
 \Ext_X^i(Q(P_2),S(P_1))\simeq H^{i-1}\bigl(A,W(aH)\otimes P_1P_2^{-1}\bigr)
\]
for all $i$, which gives $\dim\Ext_X^1(Q(P_2),S(P_1))=h$ and the vanishing for $i\ne1$ by~\eqref{eq:oddvanish}. With the arguments reversed, $\pi_*\OO_X(3)=\Sym^3\VV$ and $R^1\pi_*\OO_X(3)=0$, so the groups $\Ext_X^i(S(P_1),Q(P_2))$ are the cohomology groups of $\Sym^3\VV(-aH)\otimes P_2P_1^{-1}$, yielding the asserted vanishings.

A line-bundle extension is locally split as a sequence of modules and thus has a locally free middle term. Its Chern classes follow by multiplying the total Chern classes of $S(P_1)$ and $Q(P_2)$. Finally, a nonzero extension corresponds to a nonzero section of $W(aH)\otimes P_1P_2^{-1}$, so its restriction to the generic fiber $\PP^1$ is a nonsplit extension
\[
 0\to\OO_{\PP^1}(-1)\to E|_{\PP^1}\to\OO_{\PP^1}(2)\to0.
\]
By Grothendieck's splitting theorem, $E|_{\PP^1}\simeq\OO_{\PP^1}(\alpha)\oplus\OO_{\PP^1}(\beta)$ for integers $\alpha\le\beta$ with $\alpha+\beta=1$. Twisting by $\OO_{\PP^1}(-2)$ gives
\[
 0\to\OO_{\PP^1}(-3)\to E|_{\PP^1}(-2)\to\OO_{\PP^1}\to0.
\]
Since $H^0(\OO_{\PP^1}(-3))=0$ and the connecting map $H^0(\OO_{\PP^1})\to H^1(\OO_{\PP^1}(-3))$ is nonzero (hence injective, as $h^0(\OO_{\PP^1})=1$), we obtain $h^0(E|_{\PP^1}(-2))=0$. Thus $\beta-2<0$, i.e., $\beta\le1$; combined with $\alpha+\beta=1$ and $\alpha\le\beta$, this forces $(\alpha,\beta)=(0,1)$, so $E|_{\PP^1}\simeq\OO_{\PP^1}\oplus\OO_{\PP^1}(1)$.
\end{proof}

Lemma~\ref{lem:oddExt} determines the splitting type of $E$ on the generic fiber. We next identify the locus on $A$ where the splitting type jumps and express it via the higher direct image $R^1\pi_*(E(-1))$, whose associated pushforward sequence will also govern stability in Section~\ref{sec:nefexact}.

\begin{proposition}\label{prop:jumping}
Let $E$ be the vector bundle corresponding to a nonzero section $\xi\in H^0(A,\allowbreak W(aH)\otimes P_1P_2^{-1})$ via~\eqref{eq:oddextension}, and let $J=Z(\xi)\subset A$ be the zero scheme of $\xi$. For every closed point $x\in A$,
\[
 E|_{\pi^{-1}(x)}\simeq
 \begin{cases}
 \OO_{\PP^1}\oplus\OO_{\PP^1}(1),&x\notin J,\\
 \OO_{\PP^1}(-1)\oplus\OO_{\PP^1}(2),&x\in J.
 \end{cases}
\]
More intrinsically, setting $R=H^{\otimes a}\otimes P_1\otimes(\det\VV)^{-1}$, there is a canonical isomorphism
\[
 R^1\pi_*(E(-1))\simeq R\otimes\OO_J.
\]
\end{proposition}

\begin{proof}
Twisting~\eqref{eq:oddextension} by $\OO_X(-1)$ gives
\[
 0\to\OO_X(-2)\otimes\pi^*(H^{\otimes a}\otimes P_1)\to E(-1)\to\OO_X(1)\otimes\pi^*P_2\to0.
\]
Since $\pi_*\OO_X(-2)=0$, $R^1\pi_*\OO_X(-2)\simeq(\det\VV)^{-1}$ (as shown in the proof of Lemma~\ref{lem:oddExt}), $\pi_*\OO_X(1)=\VV$, and $R^1\pi_*\OO_X(1)=0$, the projection formula and the long exact sequence of direct images yield
\begingroup
\setlength{\abovedisplayskip}{4pt plus 1pt minus 1pt}
\setlength{\belowdisplayskip}{4pt plus 1pt minus 1pt}
\setlength{\abovedisplayshortskip}{0pt}
\setlength{\belowdisplayshortskip}{4pt plus 1pt minus 1pt}
\[
 0\to\pi_*E(-1)\to\VV\otimes P_2\xrightarrow{\partial}R\to R^1\pi_*E(-1)\to0.
\]
\endgroup
Under the pairing of Lemma~\ref{lem:oddExt}, we have
\[
 \Ext_X^1(Q(P_2),S(P_1))\simeq H^0\bigl(A,\VV^*\otimes(\det\VV)^{-1}\otimes H^{\otimes a}\otimes P_1P_2^{-1}\bigr)=\Hom_A(\VV\otimes P_2,R),
\]
and the connecting homomorphism $\partial$ is given by contraction with $\xi$. Locally on $A$, trivializing $\VV\otimes P_2$ and $R$ represents $\partial$ by the two coordinate components of $\xi$, which generate $\II_J$. Hence $\operatorname{Im}(\partial)=R\otimes\II_J$, and taking cokernels gives $R^1\pi_*(E(-1))\simeq R/(R\otimes\II_J)\simeq R\otimes\OO_J$.

Finally, restricting~\eqref{eq:oddextension} to the fiber $\pi^{-1}(x)\simeq\PP^1$ over a closed point $x\in A$ gives an extension
\[
 0\to\OO_{\PP^1}(-1)\to E|_{\pi^{-1}(x)}\to\OO_{\PP^1}(2)\to0.
\]
Since $\pi$ has one-dimensional fibers, cohomology and base change for $R^1\pi_*\OO_X(-3)$ identifies the extension class of this fiber sequence with the fiber value $\xi(x)$. Thus the sequence splits as $\OO_{\PP^1}(-1)\oplus\OO_{\PP^1}(2)$ when $x\in J$ (i.e., $\xi(x)=0$), and is nonsplit---hence isomorphic to $\OO_{\PP^1}\oplus\OO_{\PP^1}(1)$ by the splitting argument in Lemma~\ref{lem:oddExt}---when $x\notin J$.
\end{proof}

\section{Stability, the moduli component, and the upper wall}\label{sec:nefexact}

This section proves Theorem~\ref{thm:introA}. In Section~\ref{subsec:exact} we determine exactly when an extension is slope-stable and show that the common chamber is optimal. In Section~\ref{subsec:component} we show that the extensions form a smooth projective connected component of the moduli space and transfer its Brill--Noether schemes to $A$. In Section~\ref{subsec:product} we treat the product $A\times\PP^1$, where the component persists at its upper slope wall with respect to Gieseker stability.

\subsection{Exact stability}\label{subsec:exact}

To test the $\mu_{L_t}$-stability of a middle term $E$ of~\eqref{eq:oddextension}, we must bound the slopes of all saturated rank-one subsheaves of $E$. On the smooth threefold $X$, any such subsheaf is reflexive and hence a line bundle $\OO_X(j)\otimes\pi^*N$, with $j\le1$ by the general fiber splitting of Lemma~\ref{lem:oddExt}. Write $c=c_1(\VV)$. For a nonzero extension class $\xi$, let $D$ be the divisorial part of its zero scheme (allowing $D=0$) as a section of $\VV^*\otimes(\det\VV)^{-1}\otimes H^{\otimes a}\otimes P_1P_2^{-1}$. Set
\[
 \beta_t=c+2tH,\quad h_t=H\cdot\beta_t,\quad
 c_t=c\cdot\beta_t,\quad A_t=v_0+2tu+t^2d.
\]
Thus $A_t=\zeta L_t^2$ and $\pi^*N\cdot L_t^2=N\cdot\beta_t$. Since $c$ is nef and $H$ is ample, $\beta_t$ is ample for all $t>0$. We first isolate the unique maximal candidate subsheaf of relative degree one, and then rule out destabilizing subsheaves of relative degree $j\le0$ other than $S(P_1)$.

\begin{lemma}\label{lem:canonicalT}
Let $E$ be a nonzero extension in~\eqref{eq:oddextension} of class $\xi$ with divisorial zero locus $D\subset A$. Define $T_\xi:=\OO_X(1)\otimes\pi^*\pi_*E(-1)$.
\begin{enumerate}[label=(\roman*)]
\item There is an isomorphism
\begin{equation}\label{eq:nefcanonical}
 T_\xi\simeq\OO_X(1)\otimes\pi^*\bigl((\det\VV)^2\otimes H^{\otimes(-a)}
 \otimes P_1^{-1}P_2^2\otimes\OO_A(D)\bigr).
\end{equation}
\item The natural evaluation map $T_\xi\hookrightarrow E$ realizes $T_\xi$ as a saturated invertible subsheaf of $E$.
\item Every invertible subsheaf $\OO_X(1)\otimes\pi^*N\hookrightarrow E$ factors through $T_\xi$, and hence satisfies $\mu_{L_t}(\OO_X(1)\otimes\pi^*N)\le\mu_{L_t}(T_\xi)$ for every $t\in\mathbb{Q}_{>0}$.
\end{enumerate}
\end{lemma}

\begin{proof}
(i) By Proposition~\ref{prop:jumping}, pushing forward~\eqref{eq:oddextension} twisted by $\OO_X(-1)$ yields the exact sequence
\[
 0\longrightarrow\pi_*E(-1)\longrightarrow\VV\otimes P_2\xrightarrow{\ \partial\ }R:=H^{\otimes a}\otimes(\det\VV)^{-1}\otimes P_1,
\]
whose image is $R\otimes\II_J$. Decomposing $J=Z(\xi)$ into its divisorial part $D$ and zero-dimensional (or empty) part $Z$ gives $\II_J=\OO_A(-D)\II_Z$, so $\operatorname{Im}(\partial)=R(-D)\II_Z$. Locally on $A$, trivializing $\VV\otimes P_2$ and $R$ represents $\partial$ by two regular functions; dividing by their greatest common divisor (which locally defines $D$) leaves coprime elements $f,g$ of the factorial local ring $\OO_{A,x}$. Thus $\ker(\partial)$ is freely generated by $(-g,f)^{\mathsf T}$, so $\pi_*E(-1)$ is invertible. Since $Z$ has codimension two, $\det\bigl(R(-D)\II_Z\bigr)\simeq R(-D)$. Taking determinants in $0\to\pi_*E(-1)\to\VV\otimes P_2\to R(-D)\II_Z\to0$ gives
\[
 \pi_*E(-1)\simeq\det(\VV\otimes P_2)\otimes R(-D)^{-1}
 =(\det\VV)^2\otimes H^{\otimes(-a)}\otimes P_1^{-1}P_2^2\otimes\OO_A(D),
\]
and tensoring the pullback with $\OO_X(1)$ yields~\eqref{eq:nefcanonical}.

(ii) Adjunction applied to the identity of $\pi_*E(-1)$ gives the evaluation map $\pi^*\pi_*E(-1)\to E(-1)$, and twisting by $\OO_X(1)$ gives $T_\xi\hookrightarrow E$. Because $X$ is smooth and $T_\xi$ is invertible inside the locally free sheaf $E$, the quotient $E/T_\xi$ is torsion-free---that is, $T_\xi$ is saturated---if and only if $T_\xi\to E$ does not vanish identically along any prime divisor of $X$. Every prime divisor on $X$ is either horizontal or vertical:
\begin{itemize}
\item For a general point $x\notin J$, Proposition~\ref{prop:jumping} gives $R^1\pi_*E(-1)=0$ near $x$, so $\pi_*E(-1)\otimes k(x)\simeq H^0(\pi^{-1}(x),E(-1)|_{\pi^{-1}(x)})\simeq\CC$ by cohomology and base change. On the fiber $\pi^{-1}(x)\simeq\PP^1$, $T_\xi\to E$ restricts to the inclusion $\OO_{\PP^1}(1)\hookrightarrow\OO_{\PP^1}\oplus\OO_{\PP^1}(1)$, which is nonzero; hence $T_\xi\to E$ cannot vanish along any horizontal prime divisor.
\item Under adjunction, the composite $T_\xi\to E\to Q(P_2)$ corresponds to the inclusion $\pi_*E(-1)\hookrightarrow\pi_*(Q(P_2)(-1))=\VV\otimes P_2$. Locally on $A$, this map is given by the coprime column $(-g,f)^{\mathsf T}$, which vanishes only in codimension two and never along a curve $\Gamma\subset A$. Thus $T_\xi\to E$ cannot vanish along any vertical prime divisor $\pi^{-1}(\Gamma)$.
\end{itemize}
Hence $T_\xi$ is saturated.

(iii) By adjunction, any injection $\OO_X(1)\otimes\pi^*N\hookrightarrow E$ corresponds to a nonzero injection $N\hookrightarrow\pi_*E(-1)$; pulling back and twisting by $\OO_X(1)$ shows that it factors through $T_\xi$. Since $\pi_*E(-1)\otimes N^{-1}$ is effective and $\beta_t$ is ample, $(c_1(\pi_*E(-1))-c_1(N))\cdot\beta_t\ge0$, giving the slope inequality.
\end{proof}

Having isolated the unique maximal candidate subsheaf of relative degree one, we next use the nefness of $\VV$ to rule out destabilizing subsheaves of nonpositive relative degree outside $S(P_1)$.

\begin{lemma}\label{lem:nefmaxslope}
For every $t\in\mathbb{Q}_{>0}$ and every integer $q\ge0$, the maximal slope with respect to $\beta_t$ satisfies
\begin{equation}\label{eq:nefmaxslope}
 \mu_{\beta_t,\max}(\VV)\le c_t,\qquad
 \mu_{\beta_t,\max}(\Sym^q\VV)\le q c_t.
\end{equation}
Consequently, if $aH-3c$ is ample and $t\in\mathbb{Q}_{>0}$ satisfies $3A_t\ge ah_t$, then for any extension $E$ as in~\eqref{eq:oddextension}, every saturated rank-one subsheaf $T\subset E$ of relative degree $j\le0$ not contained in $S(P_1)$ satisfies $\mu_{L_t}(T)<\mu_{L_t}(E)$.
\end{lemma}

\begin{proof}
Since $c$ is nef and $\beta_t$ is ample, $c_t=c\cdot\beta_t\ge0$. If $\VV$ is $\beta_t$-semistable, then $\mu_{\beta_t,\max}(\VV)=c_t/2\le c_t$. If $\VV$ is unstable with Harder--Narasimhan filtration $0\to L_1\to\VV\to L_2\to0$, restricting the rank-one torsion-free quotient $L_2$ to a general smooth curve in a sufficiently divisible multiple of $\beta_t$ avoiding its finite singular set shows by nefness of $\VV$ that $c_1(L_2)\cdot\beta_t\ge0$; hence $\mu_{\beta_t,\max}(\VV)=c_1(L_1)\cdot\beta_t=(c-c_1(L_2))\cdot\beta_t\le c_t$. The second inequality in~\eqref{eq:nefmaxslope} follows from the characteristic-zero semistability of tensor products~\cite{HL} applied to the Harder--Narasimhan factors of $\VV$, together with $\Sym^q\VV$ being a direct summand of $\VV^{\otimes q}$.

Now write $T=\OO_X(j)\otimes\pi^*N$ with $j\le0$. If $T\not\subset S(P_1)=\ker(E\to Q(P_2))$, then the composite $T\to Q(P_2)$ is nonzero. Twisting by $T^{-1}$ and pushing forward along $\pi$ (using $\pi_*\OO_X(2-j)=\Sym^{2-j}\VV$ for $2-j\ge2$) gives a nonzero section in $H^0(A,\Sym^{2-j}\VV\otimes P_2N^{-1})$. This embeds $N\otimes P_2^{-1}$ into $\Sym^{2-j}\VV$, so~\eqref{eq:nefmaxslope} yields $N\cdot\beta_t\le(2-j)c_t$. Adding the relative contribution $jA_t$, we get
\[
 \mu_{L_t}(T)=jA_t+N\cdot\beta_t\le jA_t+(2-j)c_t=j(A_t-c_t)+2c_t.
\]
Since $aH-3c$ and $\beta_t$ are ample, $ah_t>3c_t$. Whenever $3A_t\ge ah_t$, we have $A_t>c_t$, and since $j\le0$,
\[
 \mu_{L_t}(T)\le2c_t<\frac{2}{3}ah_t\le\frac{A_t+ah_t}{2}=\mu_{L_t}(E). \qedhere
\]
\end{proof}

Combining Lemmas~\ref{lem:canonicalT} and~\ref{lem:nefmaxslope} reduces slope stability to testing only $S(P_1)$ and $T_\xi$, yielding an exact criterion and a uniform stability chamber.

\begin{theorem}\label{thm:nefexact}
Assume that $aH-3c$ is ample. For any rational $t>0$, a nonzero extension $E$ as in~\eqref{eq:oddextension} is $\mu_{L_t}$-stable if and only if
\begin{equation}\label{eq:nefexact}
 3A_t>a h_t\qquad\text{and}\qquad A_t+4c_t+2D\cdot\beta_t<3a h_t,
\end{equation}
and is $\mu_{L_t}$-semistable if and only if the corresponding non-strict inequalities hold. The two stability walls are governed respectively by the line subbundle $S(P_1)$ and by the canonical saturated invertible subsheaf $T_\xi$ of~\eqref{eq:nefcanonical}. In particular, every nonzero extension~\eqref{eq:oddextension} is $\mu_{L_t}$-stable for all $t\in\mathbb{Q}_{>0}$ in the common chamber
\begin{equation}\label{eq:nefcommon}
 \frac{a h_t}{3}<A_t<a h_t-2c_t.
\end{equation}
\end{theorem}

\begin{proof}
Let $T=\OO_X(j)\otimes\pi^*N$ be any saturated rank-one subsheaf of $E$. By Lemma~\ref{lem:oddExt}, $j\le1$. If $j=1$, then $\mu_{L_t}(T)\le\mu_{L_t}(T_\xi)$ by Lemma~\ref{lem:canonicalT}(iii). If $j\le0$ and $T\subset S(P_1)$, then $S(P_1)\otimes T^{-1}$ is effective, so $\mu_{L_t}(T)\le\mu_{L_t}(S(P_1))$. By Lemma~\ref{lem:nefmaxslope}, if $3A_t\ge ah_t$, then every $T$ with $j\le0$ and $T\not\subset S(P_1)$ satisfies $\mu_{L_t}(T)<\mu_{L_t}(E)$. Thus, in both the stable and semistable cases, $\mu_{L_t}$-stability (resp.\ $\mu_{L_t}$-semistability) is equivalent to $\max\{\mu_{L_t}(S(P_1)),\mu_{L_t}(T_\xi)\}<\mu_{L_t}(E)$ (resp.\ $\le$). Computing the two slope differences against $\mu_{L_t}(E)=(A_t+ah_t)/2$ gives
\[
 2\bigl(\mu_{L_t}(S(P_1))-\mu_{L_t}(E)\bigr)=ah_t-3A_t,
 \qquad
 2\bigl(\mu_{L_t}(T_\xi)-\mu_{L_t}(E)\bigr)=A_t+4c_t+2D\cdot\beta_t-3ah_t,
\]
which proves both equivalences in~\eqref{eq:nefexact}.

Finally, in the pushforward sequence of Lemma~\ref{lem:canonicalT}(i), $R(-D)\II_Z$ is a rank-one torsion-free quotient of the nef bundle $\VV\otimes P_2$. By the same curve-restriction argument as in Lemma~\ref{lem:nefmaxslope}, its $\beta_t$-degree is nonnegative:
\[
 (aH-c-D)\cdot\beta_t=ah_t-c_t-D\cdot\beta_t\ge0,
\]
so $D\cdot\beta_t\le ah_t-c_t$. Substituting this bound into the second inequality of~\eqref{eq:nefexact} gives
\[
 A_t+4c_t+2D\cdot\beta_t\le A_t+2c_t+2ah_t<3ah_t
\]
whenever $A_t<ah_t-2c_t$. Together with $3A_t>ah_t$, this shows that every nonzero extension~\eqref{eq:oddextension} is $\mu_{L_t}$-stable on~\eqref{eq:nefcommon}.
\end{proof}

The following corollary shows that the common chamber~\eqref{eq:nefcommon} is optimal whenever $\VV$ admits a numerically trivial quotient line bundle.

\begin{corollary}\label{cor:nefoptimal}
Under the hypotheses of Theorem~\ref{thm:nefexact}, suppose that $\VV$ admits a numerically trivial quotient line bundle $\VV\twoheadrightarrow N$. Then the interval~\eqref{eq:nefcommon} is precisely the set of rational $t>0$ for which every nonzero extension~\eqref{eq:oddextension} is $\mu_{L_t}$-stable. In particular, this optimality holds whenever $\VV=\OO_A^{\oplus2}$, $\VV=\OO_A\oplus H^{\otimes e}$ ($e\ge0$), or $\VV$ is $\widehat A$-filtered.
\end{corollary}

\begin{proof}
By Theorem~\ref{thm:nefexact}, every nonzero extension is $\mu_{L_t}$-stable on~\eqref{eq:nefcommon}, and the first inequality in~\eqref{eq:nefexact} is independent of the extension. Let $q:\VV\twoheadrightarrow N$ be a numerically trivial quotient line bundle and set $R=H^{\otimes a}\otimes(\det\VV)^{-1}\otimes P_1$. The line bundle $L':=R\otimes N^{-1}\otimes P_2^{-1}$ is numerically equivalent to $aH-c=(aH-3c)+2c$, which is ample because $aH-3c$ is ample and $c$ is nef. On the abelian surface $A$, Kodaira vanishing and Riemann--Roch give $h^0(A,L')=\chi(A,L')=\tfrac12(aH-c)^2>0$, so $L'$ admits a nonzero global section $s$.

Consider the composite homomorphism
\[
 \xi:\VV\otimes P_2\xrightarrow{q\otimes\id}N\otimes P_2\xrightarrow{\ s\ }R,
\]
which defines a nonzero extension class since $q\otimes\id$ is surjective and $s\ne0$. Because $q$ is a quotient of vector bundles, its dual $N^{-1}\hookrightarrow\VV^*$ is a subbundle (nowhere vanishing on fibers), so the zero scheme of $\xi$ coincides purely with the effective divisor $D=\operatorname{div}(s)\in|L'|$. Thus
\[
 D\cdot\beta_t=(aH-c)\cdot\beta_t=ah_t-c_t.
\]
For this extension, the second inequality in~\eqref{eq:nefexact} simplifies to $A_t<ah_t-2c_t$. Consequently, for every rational $t>0$ outside~\eqref{eq:nefcommon}, at least one nonzero extension fails to be $\mu_{L_t}$-stable.
\end{proof}

\begin{remark}\label{rem:nefwall}
The hypothesis $aH-3c$ ample is available for every fixed nef $\VV$ after increasing $a$. It holds for every $a\ge1$ when $c\equiv0$, and in the split case $c=eH$ it is precisely $a>3e$. The proof of Theorem~\ref{thm:nefexact} uses it only through the inequality $(aH-3c)\cdot\beta_t>0$, so the criterion holds for every $t$ satisfying this inequality; Corollary~\ref{cor:nefoptimal} also uses the ampleness of $aH-c$. For the nontrivial self-extension $\VV=\mathcal U_\eta$ of $\OO_A$ (see~\eqref{eq:unipotent}), it yields the individual criterion
\[
 \frac{2a}{3}<t<6a-\frac{4D\cdot H}{H^2}.
\]
For arbitrary $\VV$, the quotient of \eqref{eq:nefcanonical} has the form
\[
 E/T_\xi=\II_{\Gamma}\otimes\pi^*\bigl(H^{\otimes 2a}\otimes(\det\VV)^{-2}
 \otimes P_1^2P_2^{-1}\otimes\OO_A(-D)\bigr).
\]
Here $\Gamma$ is the zero scheme of the induced section of $E\otimes T_\xi^{-1}$; it is empty or a local complete-intersection curve. Indeed the section has no divisorial zeros, and its two local entries form a regular sequence wherever they vanish. At an upper slope wall this sequence displays two rank-one factors of equal slope. It does not by itself give a point of a Gieseker compactification: the lower Hilbert-polynomial coefficients must also be compared.
\end{remark}

\subsection{The moduli component}\label{subsec:component}

Once $t$ lies in the common chamber~\eqref{eq:nefcommon}, the relative extension spaces over $\widehat A\times\widehat A$ give a smooth projective family of $\mu_{L_t}$-stable bundles. To show that this family is an entire connected component of the moduli space, we prove that the classifying morphism is injective, by recovering $(P_1,P_2)$ and $[\xi]$ from $E$, and match the dimension of the parameter space with the Zariski tangent dimension $\dim\Ext_X^1(E,E)$.

\begin{theorem}\label{thm:oddcomponent}
Assume that $a\in\ZZ_{>0}$ satisfies the vanishing and positivity conditions~\eqref{eq:oddvanish} and~\eqref{eq:oddachoice}, and let $t\in\mathbb Q_{>0}$ lie in the common stability chamber~\eqref{eq:nefcommon} (for instance, $t=a$).
\begin{enumerate}[label=(\roman*)]
\item The middle terms of the nonzero extensions~\eqref{eq:oddextension} form a smooth projective connected component $\mathcal D_a$ of the moduli space $M_{X,L_t}(w)$, consisting entirely of rank-two vector bundles with numerical Chern classes as in Lemma~\ref{lem:oddExt}.
\item The map sending an extension $E$ to $(P_1,P_2)\in\widehat A\times\widehat A$ realizes $\mathcal D_a$ as a $\PP^{h-1}$-bundle over $\widehat A\times\widehat A$ of dimension $h+3$, and coincides with the Albanese morphism of $\mathcal D_a$.
\item Let $q_2:\mathcal D_a\to\widehat A$ denote the projection onto the second Albanese factor. Scheme-theoretically, the Brill--Noether loci on $\mathcal D_a$ are given by
\begin{equation}\label{eq:oddBNtransfer}
 \BN_X^k\cap\mathcal D_a
 =q_2^{-1}(Z_k),
\end{equation}
where $Z_k$ is the jump scheme for $h^0(A,\Sym^2\VV\otimes P)\ge k$. These loci are nonempty for $1\le k\le h^0(A,\Sym^2\VV)$; moreover, $h^0(A,\Sym^2\VV)\ge3$ whenever $\VV$ is globally generated.
\end{enumerate}
\end{theorem}

\begin{proof}
Using the normalized Poincar\'e line bundle $\mathcal P$ on $A\times\widehat A$, form the universal line bundles $\mathcal S$ and $\mathcal Q$ on $X\times(\widehat A\times\widehat A)$ restricting over $(P_1,P_2)$ to $S(P_1)$ and $Q(P_2)$. For the projection $p:X\times(\widehat A\times\widehat A)\to\widehat A\times\widehat A$, Lemma~\ref{lem:oddExt} and cohomology and base change show that $\mathcal T:=\mathcal{E}xt^1_p(\mathcal Q,\mathcal S)$ is a rank-$h$ vector bundle whose formation commutes with base change. Its one-dimensional subspaces are parametrized by the smooth projective variety
\[
 Z:=\PP_{\widehat A\times\widehat A}(\mathcal T^*),\qquad \dim Z=4+(h-1)=h+3.
\]
Pulling back $\OO_Z(-1)\hookrightarrow\mathcal T|_Z$ and twisting the kernel by $\OO_Z(1)$ yields a universal extension on $X\times Z$~\cite{Lange}, whose middle term $\mathcal E$ is locally free of rank two and hence flat over $Z$. By~\eqref{eq:oddachoice}, $aH-3c_1(\VV)$ is ample and $t=a$ satisfies~\eqref{eq:nefcommon} (as $3A_a-ah_a=a^2d+5au+3v_0$ and $ah_a-2c_a-A_a=a^2d-5au-2\int_Ac_1(\VV)^2-v_0$). Every fiber of $\mathcal E$ is therefore $\mu_{L_t}$-stable by Theorem~\ref{thm:nefexact}, inducing a proper classifying morphism $f:Z\to M_{X,L_t}(w)$.

To see that $f$ is injective on geometric points, test $E$ against $S(P_1')$ for any $P_1'\in\widehat A$. Since $\Hom_X(S(P_1'),Q(P_2))=0$ by Lemma~\ref{lem:oddExt}, applying $\Hom_X(S(P_1'),-)$ to~\eqref{eq:oddextension} gives
\[
 \Hom_X(S(P_1'),E)\simeq\Hom_X(S(P_1'),S(P_1))\simeq H^0(A,P_1P_1'^{-1}),
\]
which is nonzero if and only if $P_1'=P_1$, and then is one-dimensional. Hence $S(P_1)\subset E$ is the unique subsheaf of the form $S(P_1')$. Any isomorphism $\phi:E\xrightarrow{\sim}E'$ therefore forces $P_1'=P_1$, preserves $S(P_1)$, and induces an isomorphism $Q(P_2)\xrightarrow{\sim}Q(P_2')$; twisting by $\OO_X(-2)$ and pushing forward yields $P_2\simeq P_2'$. Moreover, $\phi$ acts on the line bundles $S(P_1)$ and $Q(P_2)$ by scalars $\lambda_1,\lambda_2\in\CC^*$, scaling the extension class by $\lambda_1\lambda_2^{-1}$ and defining the same point of $Z$.

Next, we compute the Zariski tangent space $\Ext_X^1(E,E)$. Write $S=S(P_1)$ and $Q=Q(P_2)$, and let $\mathcal A\subset\cEnd(E)$ be the subsheaf of endomorphisms preserving $S$. Since~\eqref{eq:oddextension} splits locally on $X$, $\mathcal A$ is the upper-triangular subbundle of $\cEnd(E)$ in local splittings, fitting into two exact sequences
\begin{align*}
 0&\to\cHom(Q,S)\to\mathcal A\to\OO_X^{\oplus2}\to0,\\
 0&\to\mathcal A\to\cEnd(E)\to\cHom(S,Q)\to0.
\end{align*}
In the first sequence, the connecting map $H^0(\OO_X)^{\oplus2}\to\Ext^1_X(Q,S)$ sends $(b,c)$ to $(b-c)\xi$, which has rank one as $\xi\ne0$. Since $\Ext^2_X(Q,S)=0$ by Lemma~\ref{lem:oddExt}, $H^1(\mathcal A)\to H^1(\OO_X)^{\oplus2}$ is surjective with kernel $\Ext^1_X(Q,S)/\CC\xi$ of dimension $h-1$. Using $h^1(\OO_X)=h^1(\OO_A)=2$, we find $h^1(\mathcal A)=(h-1)+4=h+3$. In the second sequence, $\Ext^i_X(S,Q)=0$ for $i=0,1$ by Lemma~\ref{lem:oddExt}, so $H^1(\mathcal A)\xrightarrow{\sim}H^1(\cEnd(E))$ and
\begin{equation}\label{eq:oddTang}
 \dim\Ext_X^1(E,E)=h+3.
\end{equation}

We now identify $f(Z)$ with a connected component. Since $f$ is proper and injective on geometric points, $f(Z)\subset M_{X,L_t}(w)$ is closed and irreducible of dimension $h+3$. At every closed point $[E]\in f(Z)$, the local dimension of $M_{X,L_t}(w)$ is at least $h+3$ while its Zariski tangent dimension is $h+3$ by~\eqref{eq:oddTang}. Thus $\OO_{M_{X,L_t}(w),[E]}$ is regular of dimension $h+3$, so a unique irreducible component passes through $[E]$ and coincides locally with $f(Z)$. Hence $\mathcal D_a:=f(Z)$ is open and closed in $M_{X,L_t}(w)$; being projective, it is also closed in the full Gieseker moduli space, and is therefore a smooth projective connected component. As $f:Z\to\mathcal D_a$ is proper and bijective in characteristic zero, it is finite and birational onto a smooth (hence normal) variety, so $f$ is an isomorphism by Zariski's Main Theorem. Because every morphism from a projective space $\PP^{h-1}$ to an abelian variety is constant, the $\PP^{h-1}$-bundle projection $\mathcal D_a\to\widehat A\times\widehat A$ is the Albanese morphism of $\mathcal D_a$ (with a base point chosen over $(0,0)$).

Finally, since $S(P_1)$ has relative degree $-1$, we have $R\pi_*S(P_1)=0$. Pushing forward the universal extension on $X\times Z$ along $\pi\times\id_Z$ yields
\[
 R(\pi\times\id_Z)_*\mathcal E\simeq(\id_A\times q_2)^*\bigl(\pr_A^*\Sym^2\VV\otimes\mathcal P\bigr)\otimes\pr_Z^*\mathcal L
\]
for a line bundle $\mathcal L$ on $Z$, where $\pr_A:A\times\widehat A\to A$ and $\pr_Z:A\times Z\to Z$ are the projections. Since $q_2:Z\to\widehat A$ is flat, flat base change and the invariance of degree-zero jump ideals under twisting by base line bundles (Section~\ref{sec:conventions}) yield $\BN_X^k\cap\mathcal D_a=q_2^{-1}(Z_k)$ scheme-theoretically. Over $P_2=\OO_A$, we have $h^0(X,E)=h^0(A,\Sym^2\VV)$, giving nonemptiness for $1\le k\le h^0(A,\Sym^2\VV)$; if $\VV$ is globally generated, so is the rank-three bundle $\Sym^2\VV$, which thus has at least three global sections.
\end{proof}

\begin{remark}\label{rem:oddalb}
The two endomorphism sequences in the proof of Theorem~\ref{thm:oddcomponent} also yield $\Ext_X^3(E,E)=0$ and
\[
 \dim\Ext_X^2(E,E)=2+h^2\bigl(A,\Sym^3\VV(-aH)\otimes P_2P_1^{-1}\bigr),
\]
which is strictly positive in general (for instance, it equals $2d+2$ on $A\times\PP^1$ with $a=1$). The nonvanishing of $\Ext_X^2(E,E)$ does not imply that actual obstructions occur: the smoothness proved in Theorem~\ref{thm:oddcomponent} shows that the deformations of $E$ are unobstructed. Here smoothness follows from matching the family dimension with $\dim\Ext_X^1(E,E)$, rather than from the vanishing of the obstruction space. Note also that by Theorem~\ref{thm:oddcomponent}(ii)--(iii), every Brill--Noether locus on $\mathcal D_a$ is pulled back from the second factor of $\operatorname{Alb}(\mathcal D_a)\simeq\widehat A\times\widehat A$, reducing its study to the cohomology jump schemes of $\Sym^2\VV$ on $A$. When $-K_X$ is nef, these schemes are described by Corollary~\ref{cor:trichotomy}.
\end{remark}

\subsection{The product \texorpdfstring{$A\times\PP^1$}{A x P1} and its upper wall}\label{subsec:product}

For the trivial $\PP^1$-bundle $X=A\times\PP^1$ and $a=1$, the common chamber is the explicit interval $2/3<t<2$. At its upper end $t=2$, the extensions given by proportional pairs of sections become strictly slope-semistable. Comparing the linear terms of the reduced Hilbert polynomials, we show that every member of the component nevertheless remains Gieseker-stable there.

Take $\VV=\OO_A^{\oplus2}$ and $a=1$, so $X=A\times\PP^1$. Then \eqref{eq:oddextension} becomes
\begin{equation}\label{eq:oddproductextension}
 0\to\pi^*(H\otimes P_1)\otimes\OO_X(-1)\to E\to\pi^*P_2\otimes\OO_X(2)\to0,
\end{equation}
with numerical Chern classes $(\zeta+B,2\zeta B,0)$. For these data $h=d$ by Riemann--Roch, and the vanishings in~\eqref{eq:oddvanish} follow from the positive and negative ample twists on $A$. Theorem~\ref{thm:oddcomponent} and Corollary~\ref{cor:nefoptimal} give a smooth component $\mathcal D_1$ of dimension $d+3$ throughout the optimal common chamber $2/3<t<2$.

Serre duality on $\PP^1$ identifies the extension class $\xi$ with a pair of sections $s_0,s_1$ of $M=H\otimes P_1P_2^{-1}$. Let $D$ be their maximal common effective divisor; if one section vanishes, use the divisor of the other. In the notation of Theorem~\ref{thm:nefexact}, $A_t=t^2d$, $h_t=2td$, $c_t=0$ and $D\cdot\beta_t=2t\,D\cdot H$. Hence $E$ is $\mu_{L_t}$-stable if and only if
\begin{equation}\label{eq:productchamber}
 \frac23<t<6-\frac{4D\cdot H}{d},
\end{equation}
and \eqref{eq:nefcanonical} becomes
\begin{equation}\label{eq:canonicalSubline}
 T_\xi=\OO_X(1)\otimes\pi^*\bigl(H^{-1}P_1^{-1}P_2^2\otimes\OO_A(D)\bigr).
\end{equation}
The upper bound equals $2$ exactly for proportional pairs, where $\OO_A(D)=M$ and~\eqref{eq:canonicalSubline} simplifies to $T_\xi=\OO_X(1)\otimes\pi^*P_2$. At $L_2=\zeta+2B$, we compare the reduced Hilbert polynomials $p_G(m)=\chi(G(mL_2))/\rk G$.

\begin{proposition}\label{prop:giesekerwall}
At the upper wall $L_2=\zeta+2B$ on $X=A\times\PP^1$, every nonzero extension~\eqref{eq:oddproductextension} is Gieseker-stable, so $\mathcal D_1$ persists as a smooth projective connected component of $M_{X,L_2}(w)$. In particular, when $(s_0,s_1)$ are proportional, $E$ is strictly $\mu_{L_2}$-semistable with
\begin{equation}\label{eq:wallpoly}
 p_{T_\xi}(m)-p_E(m)=-\frac d4m,
\end{equation}
whereas $T_\xi\oplus(E/T_\xi)$ is Gieseker-unstable.
\end{proposition}

\begin{proof}
Since $M=H\otimes P_1P_2^{-1}$ is numerically equivalent to $H$, we have $D\cdot H\le M\cdot H=d$, so~\eqref{eq:productchamber} shows that every nonzero extension~\eqref{eq:oddproductextension} is $\mu_{L_2}$-semistable. Equality $D\cdot H=d$ holds if and only if the residual line bundle $M(-D)$, which admits a nonzero section and satisfies $M(-D)\cdot H=0$, is trivial---that is, $\OO_A(D)\simeq M$, meaning that $s_0$ and $s_1$ vanish along the same divisor and are therefore scalar multiples of each other. In this case, substituting $\OO_A(D)=H\otimes P_1P_2^{-1}$ into~\eqref{eq:canonicalSubline} gives $T_\xi=\OO_X(1)\otimes\pi^*P_2$. By Lemma~\ref{lem:nefmaxslope} (as $3A_2=12d>4d=ah_2$), no saturated rank-one subsheaf of relative degree $j\le0$ has slope equal to $\mu_{L_2}(E)=4d$. By Lemma~\ref{lem:canonicalT}(iii), every relative-degree-one invertible subsheaf $\OO_X(1)\otimes\pi^*N\subset E$ satisfies $N\subset\pi_*E(-1)=P_2$, so equality of slopes forces $N=P_2$ and hence the subsheaf equals $T_\xi$. Thus $T_\xi$ is the unique saturated rank-one subsheaf of $E$ of equal slope.

On $X=A\times\PP^1$, we have $mL_2=m\zeta+2mB=\OO_{\PP^1}(m)\boxtimes\OO_A(2mH)$. Using $\chi(\PP^1,\OO_{\PP^1}(r))=r+1$ and Riemann--Roch on the abelian surface $A$ (namely $\chi(A,\OO_A(nH)\otimes P)=\tfrac d2 n^2$ for every $P\in\widehat A$), the K\"unneth formula and additivity along~\eqref{eq:oddproductextension} give
\begin{align*}
 \chi(T_\xi(mL_2))&=\chi(\PP^1,\OO_{\PP^1}(m+1))\,\chi(A,\OO_A(2mH)\otimes P_2)=\tfrac d2(m+2)(2m)^2,\\
 \chi(E(mL_2))&=\chi(S(P_1)(mL_2))+\chi(Q(P_2)(mL_2))=\tfrac d2\bigl(m(2m+1)^2+(m+3)(2m)^2\bigr).
\end{align*}
Dividing $\chi(E(mL_2))$ by $\rk E=2$ and subtracting cancels the cubic and quadratic terms and leaves
\[
 p_{T_\xi}(m)-p_E(m)=\frac d4\Bigl(2(m+2)(2m)^2-m(2m+1)^2-(m+3)(2m)^2\Bigr)=-\frac d4m<0
\]
for all $m>0$, proving~\eqref{eq:wallpoly}. Every rank-one subsheaf $F\subset E$ embeds into its saturation $\overline F\subset E$, which is a line bundle on the smooth threefold $X$, and $p_F(m)\le p_{\overline F}(m)$ for $m\gg0$. If $\overline F=T_\xi$, then~\eqref{eq:wallpoly} gives $p_F(m)\le p_{T_\xi}(m)<p_E(m)$ for $m\gg0$. If $\overline F\ne T_\xi$, the uniqueness of the equal-slope saturated subsheaf gives $\mu_{L_2}(\overline F)<\mu_{L_2}(E)$; hence $p_{\overline F}(m)$ has strictly smaller quadratic coefficient than $p_E(m)$, and again $p_F(m)<p_E(m)$ for $m\gg0$. Thus $E$ is Gieseker-stable. When $(s_0,s_1)$ are not proportional, $E$ is already $\mu_{L_2}$-stable. The parameter recovery and tangent-space calculation of Theorem~\ref{thm:oddcomponent} now show that $\mathcal D_1$ remains a smooth projective connected component of $M_{X,L_2}(w)$, and the Brill--Noether transfer~\eqref{eq:oddBNtransfer} is independent of the polarization.

Finally, setting $Q'=E/T_\xi$, additivity of Euler characteristics gives $p_{Q'}(m)-p_E(m)=\tfrac d4m>0$ for $m>0$, so $Q'$ Gieseker-destabilizes $T_\xi\oplus Q'$.
\end{proof}

The same K\"unneth calculation applies for every integer $a\ge1$: on $X=A\times\PP^1$, the common chamber is $\tfrac{2a}{3}<t<2a$, and at the upper wall $L_{2a}=\zeta+2aB$ the proportional pairs still satisfy $p_{T_\xi}(m)-p_E(m)=-\tfrac{a^2d}{4}m<0$, so $\mathcal D_a$ persists as a smooth projective connected component of $M_{X,L_{2a}}(w)$.

\section{Cohomology jump schemes on abelian varieties}\label{sec:jump}

This section works on an abelian variety alone and proves Theorem~\ref{thm:introB}. Since the argument does not use the dimension, we let $A$ be a complex abelian variety of dimension $g$ throughout the section; in the rest of the paper $g=2$. We consider families $\{F\otimes P\}_{P\in\widehat A}$, where $F$ is a nef semihomogeneous bundle, and determine their cohomology jump schemes. Section~\ref{subsec:homogjump} treats homogeneous bundles through the finite-length Fourier--Mukai transform. Section~\ref{subsec:twistjump} treats homogeneous bundles twisted by a nef line bundle $L$: the jump schemes are pulled back, along the restriction map $\widehat A\to\widehat K$ to the dual of the abelian subvariety $K=K(L)^0$, from Fitting schemes of a finite-length module. Section~\ref{subsec:semihom} reduces arbitrary nef semihomogeneous bundles to this case by an isogeny, which proves Theorem~\ref{thm:introB}(i). Section~\ref{subsec:surfacejump} computes the multiplicities along elliptic directions and proves Theorem~\ref{thm:introB}(ii) for abelian surfaces.

\subsection{Homogeneous bundles}\label{subsec:homogjump}

We first recall how the Fourier--Mukai transform turns homogeneous vector bundles on $A$ into finite-length sheaves on $\widehat A$, and use this to express their cohomology jump ideals as Fitting ideals of local $\Ext$-modules.

Call a vector bundle $\mathcal W$ on $A$ \emph{$\widehat A$-filtered} if it has a filtration whose successive quotients are line bundles in $\widehat A$. Using the normalized Poincar\'e bundle $\mathcal P$ on $A\times\widehat A$ (with the corresponding definition in arbitrary dimension $g$), write $P_y=\mathcal P|_{A\times\{y\}}$, let $p:A\times\widehat A\to\widehat A$ be the projection, and set $\Phi_A(F)=Rp_*(\pr_A^*F\otimes\mathcal P)$. We use the convention $H^i(\mathcal C[n])=H^{i+n}(\mathcal C)$. For a closed point $y\in\widehat A$ write $R_y=\OO_{\widehat A,y}$, a regular local ring of dimension $g$.
By Miyanishi and Mukai, $\widehat A$-filtered bundles are exactly the homogeneous, that is translation-invariant, vector bundles~\cite{Miyanishi,Mukai78}. Their Fourier--Mukai transforms give the equivalence with finite-length sheaves on $\widehat A$~\cite[Example~3.2]{Mukai81}. Such bundles are nef.

\begin{proposition}\label{prop:FMfitting}
Let $\mathcal W$ be an $\widehat A$-filtered vector bundle of rank $r$ on $A$. Then
\[
 \Phi_A(\mathcal W)\simeq M_{\mathcal W}[-g],
\]
where $M_{\mathcal W}$ is a coherent sheaf of finite length $r$ on $\widehat A$, supported at those points $y\in\widehat A$ for which $P_y\simeq N^{-1}$ for some graded piece $N$ of $\mathcal W$. For every $y\in\widehat A$ and $k\ge1$, the degree-zero cohomology jump ideal of the family $\{\mathcal W\otimes P_{y'}\}_{y'\in\widehat A}$ defining the locus $\{h^0\ge k\}$ is given in $R_y$ by the Fitting ideal
\[
 \operatorname{Fitt}_{k-1}\bigl(\Ext^g_{R_y}(M_{\mathcal W,y},R_y)\bigr).
\]
In particular, $h^0(A,\mathcal W\otimes P_y)$ equals the dimension of the socle of the $R_y$-module $M_{\mathcal W,y}$.
\end{proposition}

\begin{proof}
For $N\in\widehat A$, let $y_N\in\widehat A$ be the point with $P_{y_N}\simeq N^{-1}$, and let $t(y)=y-y_N$ be translation by $-y_N$ on $\widehat A$. By the theorem of the square, $\pr_A^*N\otimes\mathcal P\simeq(1\times t)^*\mathcal P$. Flat base change along $p$ and $1\times t$ together with Mukai's identity $Rp_*\mathcal P\simeq k(0)[-g]$~\cite[proof of Theorem~2.2]{Mukai81} then yields
\[
 \Phi_A(N)\simeq t^*(Rp_*\mathcal P)\simeq t^*\bigl(k(0)[-g]\bigr)=k(y_N)[-g].
\]
Given an extension $0\to\mathcal W'\to\mathcal W\to N\to0$ with $N\in\widehat A$, tensoring its pullback to $A\times\widehat A$ with $\mathcal P$ and taking higher direct images along $p$ yields a long exact sequence of sheaves on $\widehat A$. Since $R^ip_*(\pr_A^*N\otimes\mathcal P)=0$ for $i\ne g$ and $R^qp_*=0$ for all $q>g$ (as $p$ has $g$-dimensional fibers), induction on the filtration rank shows that $R^ip_*(\pr_A^*\mathcal W\otimes\mathcal P)=0$ for $i\ne g$ and produces exact sequences
\[
 0\longrightarrow M_{\mathcal W'}\longrightarrow M_{\mathcal W}\longrightarrow k(y_N)\longrightarrow 0.
\]
Since $k(y_N)$ has length $1$ and support $\{y_N\}$, induction on the filtration implies that $M_{\mathcal W}$ has finite length $r$ and support $\{[N^{-1}]\}$ as $N$ ranges over the graded pieces of $\mathcal W$.

Outside $\operatorname{Supp}(M_{\mathcal W})$ the complex $\Phi_A(\mathcal W)$ is acyclic and all positive jump ideals are the unit ideal. Fix a closed point $y\in\operatorname{Supp}(M_{\mathcal W})$. Since $M_{\mathcal W,y}$ is a nonzero finite-length $R_y$-module, it is annihilated by a power of $\mathfrak m_y$ and hence has depth $0$. By the Auslander--Buchsbaum formula over the $g$-dimensional regular local ring $R_y$, a minimal free resolution of $M_{\mathcal W,y}$ has length $g$:
\[
 0\longrightarrow F_g\xrightarrow{\ \varphi\ }F_{g-1}\longrightarrow\cdots\longrightarrow F_0\longrightarrow M_{\mathcal W,y}\longrightarrow 0.
\]
Placed in degrees $0,\dots,g$ with $F_g$ in degree zero, this complex represents $\Phi_A(\mathcal W)$ near $y$, and degree-zero jump ideals do not depend on the perfect representative (Section~\ref{sec:conventions}). The ideal for $h^0\ge k$ is therefore $I_{\rk F_g-k+1}(\varphi)$. Dualizing into $R_y$ yields
\[
 0\longrightarrow F_0^\vee\longrightarrow\cdots\longrightarrow F_{g-1}^\vee\xrightarrow{\ \varphi^\vee\ }F_g^\vee\longrightarrow 0.
\]
Since $R_y$ is regular (hence Cohen--Macaulay) of dimension $g$ and $\operatorname{Ann}_{R_y}(M_{\mathcal W,y})$ has height $g$, we have $\Ext^i_{R_y}(M_{\mathcal W,y},R_y)=0$ for $i<g$ and $\Ext^g_{R_y}(M_{\mathcal W,y},R_y)\simeq\operatorname{coker}(\varphi^\vee:F_{g-1}^\vee\to F_g^\vee)$. Because this module is presented with $\rk F_g^\vee=\rk F_g$ generators and relation matrix $\varphi^\vee$, its $(k-1)$-st Fitting ideal is $I_{\rk F_g-k+1}(\varphi^\vee)=I_{\rk F_g-k+1}(\varphi)$.

Finally, by derived base change to $y$, the fiber cohomology is computed by $(F_\bullet[-g])\otimes_{R_y}k(y)$. By minimality of the resolution, $\varphi\otimes k(y)=0$, so
\[
 h^0(A,\mathcal W\otimes P_y)=\rk F_g=\dim_{k(y)}\operatorname{Tor}_g^{R_y}\bigl(M_{\mathcal W,y},k(y)\bigr).
\]
On the other hand, the Koszul complex on a regular system of parameters $(x_1,\dots,x_g)$ of $R_y$ is a free resolution of $k(y)$; tensoring it with $M_{\mathcal W,y}$ identifies $\operatorname{Tor}_g^{R_y}(M_{\mathcal W,y},k(y))$ with the kernel of $M_{\mathcal W,y}\xrightarrow{(x_1,\dots,x_g)}M_{\mathcal W,y}^{\oplus g}$, which is $\operatorname{soc}_{R_y}(M_{\mathcal W,y})$.
\end{proof}

\begin{remark}\label{rem:localFitting}
The second paragraph of the proof of Proposition~\ref{prop:FMfitting} applies verbatim to a finite-length module $N$ placed in degree $s$ over a regular local ring $R$ of dimension $s$: the degree-zero jump ideals of $N[-s]$ are the Fitting ideals $\operatorname{Fitt}_{k-1}(\Ext^s_R(N,R))$, and the rank in degree zero of a minimal free representative equals $\dim\operatorname{soc}_R(N)$. We use this below with $s=\dim K$.
\end{remark}

\subsection{Twists by a nef line bundle}\label{subsec:twistjump}

Let $L$ be a nef line bundle on $A$. Put
\[
 \phi_L(a)=t_a^*L\otimes L^{-1},\qquad K=(\ker\phi_L)^0,\qquad s=\dim K,\qquad f:A\to\bar A=A/K,
\]
and let $i:K\hookrightarrow A$ be the inclusion. By~\cite[Corollary~5.4]{Gulbrandsen}, there are $P_0\in\widehat A$ and an ample line bundle $M$ on $\bar A$ with
\begin{equation}\label{eq:nefdescent}
 L\simeq f^*M\otimes P_0 .
\end{equation}
Recall from the Appell--Humbert theorem~\cite[Chapter~2]{BiL} that, writing $A=V/\Lambda$ with $V=T_0A$, the first Chern class of any line bundle $L$ on $A$ is uniquely represented by a Hermitian form $H_L$ on $V$ whose imaginary part $E_L=\operatorname{Im}H_L$ is integer-valued on $\Lambda$; extending $\QQ$-linearly attaches a Hermitian form $H_\delta$ to every class $\delta\in\NS(A)_\QQ$. The class $\delta$ is nef if and only if $H_\delta$ is positive semidefinite, and the derivative of $\phi_L:A\to\widehat A$ at $0$ has kernel $\ker H_L$, so $T_0K=\ker H_L$~\cite[Section~2.4]{BiL}. In~\eqref{eq:nefdescent}, $H_L$ descends to a positive definite Hermitian form on $V/T_0K=T_0\bar A$ with integral imaginary part, yielding the ample line bundle $M$, while $L\otimes f^*M^{-1}\in\widehat A$ has trivial first Chern class. The factor $P_0$ cannot be omitted, since a nontrivial element of $\widehat A$ is nef. Put $m=h^0(\bar A,M)=\chi(\bar A,M)\ge1$; if $s=g$, then $\bar A$ is a point and $m=1$. Let $\rho=\widehat\imath:\widehat A\to\widehat K$ be the restriction map. It is smooth and surjective, with kernel the abelian subvariety $\widehat f(\widehat{\bar A})$ of dimension $g-s$. For $y\in\widehat K$ write $\tau_y=\rho^{-1}(y)$ and $R_y=\widehat\OO_{\widehat K,y}$.

We use three identities for the normalized Poincar\'e transforms:
\begin{gather}
 \Phi_A(F\otimes^{\mathbf L}G)\simeq\Phi_A(F)\star\Phi_A(G)\,[g],\label{eq:FMconv}\\
 \Phi_D\circ Lh^*\simeq R\widehat h_*\circ\Phi_{D'}\,[\dim D'-\dim D],\label{eq:FMhom}\\
 \Phi_A(F\otimes P_a)\simeq t_a^*\Phi_A(F),\qquad t_a(\alpha)=\alpha+a.\label{eq:FMtrans}
\end{gather}
Here $\star$ is the derived Pontryagin product, the pushforward of the external product along the addition map, and $h:D\to D'$ is a homomorphism of abelian varieties. Identities~\eqref{eq:FMconv} and~\eqref{eq:FMtrans} are~\cite[(3.7) and (3.1)]{Mukai81}. For~\eqref{eq:FMhom}, the Poincar\'e kernels give $\Phi_{D'}Rh_*\simeq L\widehat h^*\Phi_D$ by flat base change and the projection formula; applying this to $\widehat h$ and using Mukai's inversion on both sides gives~\eqref{eq:FMhom}; see also~\cite[proof of Proposition~5.1]{Gulbrandsen} and~\cite[Section~11.3]{Polishchuk}. The shifts are fixed by the case of structure sheaves. Note that $\Phi_A(P_a)=k(-a)[-g]$, so $t_a^*$ moves a support point $z$ to $z-a$.

\begin{theorem}\label{thm:twistjump}
Let $\mathcal W$ be a homogeneous bundle of rank $r$ on $A$, let $L$ be a nef line bundle, and keep the notation above. Write $Z_k$ for the degree-zero jump scheme of the family $\{\mathcal W\otimes L\otimes P\}_{P\in\widehat A}$ for $h^0\ge k$, and put
\[
 N=M_{(\mathcal W\otimes L)|_K},\qquad N_y^\vee=\Ext^s_{R_y}(N_y,R_y),
\]
a finite-length sheaf on $\widehat K$ and its local duals.
\begin{enumerate}[label=(\roman*)]
\item $\Phi_A(\mathcal W\otimes L)\simeq\mathcal T[-s]$ for a coherent sheaf $\mathcal T$ supported on finitely many translates $\tau_y$, and the completion of $\mathcal T$ at every point of $\tau_y$ is isomorphic to $(\widehat\OO_{\widehat A}\otimes_{R_y}N_y)^{\oplus m}$. In particular $H^i(\mathcal W\otimes L\otimes P)=0$ for $i>s$ and all $P$.
\item Along $\tau_y$, the ideal of $Z_k$ is the pullback under $\rho$ of $\operatorname{Fitt}_{k-1}\bigl((N_y^\vee)^{\oplus m}\bigr)\subset R_y$. Thus $Z_k$ is the preimage under $\rho$ of a finite subscheme of $\widehat K$, and every nonempty $Z_k\ne\widehat A$ has codimension exactly $s$.
\item For every $\alpha\in\tau_y$, $h^0(A,\mathcal W\otimes L\otimes P_\alpha)=m\dim_\CC\operatorname{soc}_{R_y}(N_y)$.
\end{enumerate}
In particular, if $L$ is numerically trivial, then $s=g$, $\rho$ is the identity, $m=1$, and the ideals of $Z_k$ are those of Proposition~\ref{prop:FMfitting} applied to $\mathcal W\otimes L$. If $L$ is ample, then $s=0$ and $Z_k=\widehat A$ for $k\le rm$, with $m=\chi(L)$, and $Z_k=\varnothing$ otherwise.
\end{theorem}

\begin{proof}
Put $M'=M_{\mathcal W\otimes P_0}$, which is a translate of $M_{\mathcal W}$ by~\eqref{eq:FMtrans}. Since $M$ is ample on $\bar A$, Kodaira vanishing gives $H^i(\bar A,M\otimes P)=0$ for $i>0$ and $h^0(\bar A,M\otimes P)=\chi(\bar A,M)=m$ for all $P\in\widehat{\bar A}$, so $V_M:=\Phi_{\bar A}(M)$ is a vector bundle of rank $m$ on $\widehat{\bar A}$ in degree zero. Since $f:A\to\bar A$ has connected kernel $K$, the dual homomorphism $\widehat f:\widehat{\bar A}\hookrightarrow\widehat A$ is a closed immersion, so~\eqref{eq:FMhom} for $f$ gives $\Phi_A(f^*M)\simeq\widehat f_*V_M[-s]$. By~\eqref{eq:nefdescent} and~\eqref{eq:FMconv},
\[
 \Phi_A(\mathcal W\otimes L)\simeq M'[-g]\star\widehat f_*V_M[-s]\,[g]=\mathcal T[-s],\qquad \mathcal T:=\sigma_*(M'\boxtimes V_M),
\]
where $\sigma(z,b)=z+\widehat f(b)$. Because $\operatorname{Supp}(M')$ is finite and $\widehat f$ is a closed immersion, $\sigma$ is finite on $\operatorname{Supp}(M')\times\widehat{\bar A}$, so $R^i\sigma_*=0$ for $i>0$ and $\mathcal T$ is supported on the translates $\tau_{\rho(z)}$, $z\in\operatorname{Supp}(M')$. This gives the vanishing in~(i).

To relate $\mathcal T$ to the restriction $(\mathcal W\otimes L)|_K$, we apply~\eqref{eq:FMhom} to the inclusion $i:K\hookrightarrow A$, which yields $\Phi_K(i^*F)\simeq R\rho_*\Phi_A(F)[g-s]$. Taking $F=\mathcal W\otimes P_0$, since $f^*M|_K$ is trivial and $M'$ is zero-dimensional (so $R^q\rho_*M'=0$ for $q>0$), we obtain
\[
 N=M_{(\mathcal W\otimes L)|_K}=M_{(\mathcal W\otimes P_0)|_K}\simeq\rho_*M'.
\]
Hence, for $y\in\widehat K$ and $S_y=\{z\in\operatorname{Supp}(M'):\rho(z)=y\}$,
\begin{equation}\label{eq:Nsum}
 N_y\simeq\bigoplus_{z\in S_y}\operatorname{Res}_{R_y}M'_z,
\end{equation}
where $\operatorname{Res}_{R_y}$ denotes restriction of scalars along $\rho$.

We next determine the formal completion of $\mathcal T$ at an arbitrary point $\alpha\in\tau_y$. In translated formal logarithm coordinates at $\alpha$, choose parameters $y=(y_1,\dots,y_s)$ pulled back from $\widehat K$ (so $R_y\simeq\CC[[y]]$) and coordinates $u=(u_1,\dots,u_{g-s})$ restricting to formal parameters along the fiber $\tau_y$, giving $\widehat\OO_{\widehat A,\alpha}\simeq\CC[[u,y]]$. For each $z\in S_y$, the relation $\sigma(z,b)=\alpha$ determines a unique $b\in\widehat{\bar A}$ with $\widehat f(b)=\alpha-z$. Choose compatible translated formal logarithm coordinates $(u_z,y_z)$ at $z$ and $v$ at $b$, and trivialize $V_M\simeq\OO^{\oplus m}$ near $b$. In these coordinates the pullback along $\sigma$ is
\[
 \sigma^*(u)=u_z+v,\qquad \sigma^*(y)=y_z.
\]
Let $U=(U_1,\dots,U_{g-s})$ and $Y=(Y_1,\dots,Y_s)$ denote the commuting nilpotent endomorphisms of $M'_z$ induced by $u_z$ and $y_z$. On $M'_z[[v]]$, the resulting $\CC[[u,y]]$-module structure has $u_i$ acting as $v_i+U_i$ and $y_j$ acting as $Y_j$. The substitution
\[
 \Theta_z:M'_z[[v]]\longrightarrow M'_z[[u]],\qquad
 \sum_\nu v^\nu a_\nu\longmapsto\sum_\nu(u-U)^\nu a_\nu
\]
is well-defined because the ideal generated by the $U_i$ acts nilpotently. Its inverse is the substitution $u\mapsto v+U$. It is a $\CC[[u,y]]$-linear isomorphism when $u_i$ acts on the target by ordinary multiplication and $y_j$ acts by $Y_j$: indeed, $\Theta_z((v_i+U_i)w)=u_i\Theta_z(w)$, and $\Theta_z$ commutes with every $Y_j$. Thus each of the $m$ summands supplied by the trivialization of $V_M$ is isomorphic to $\CC[[u]]\widehat\otimes_\CC M'_z$. Summing over $z\in S_y$ and using~\eqref{eq:Nsum},
\[
 \widehat{\mathcal T}_\alpha\simeq\bigl(\CC[[u]]\widehat\otimes_\CC N_y\bigr)^{\oplus m}=\bigl(\widehat\OO_{\widehat A,\alpha}\otimes_{R_y}N_y\bigr)^{\oplus m}.
\]
Distinct support points are direct summands of $M'$, so they remain direct summands after pushforward even when their images coincide. This proves~(i).

It remains to compute the cohomology jump ideals and identify the resulting global schemes. Define the coherent sheaf and closed subscheme on $\widehat K$ by
\[
 \mathcal N^\vee:=\mathcal{E}\!xt^s_{\OO_{\widehat K}}(N,\OO_{\widehat K}),\qquad
 \mathcal I_{D_k}:=\operatorname{Fitt}_{k-1}\bigl((\mathcal N^\vee)^{\oplus m}\bigr).
\]
The sheaf $\mathcal N^\vee$ has finite support, so $D_k$ is a finite subscheme, possibly empty. Since $\widehat K$ is smooth, $N$ locally admits a finite resolution by finite free modules; flat base change to the completion therefore identifies the completed stalk of $\mathcal N^\vee$ at $y$ with $N_y^\vee$. Fitting ideals commute with base change, so the completed ideal of $D_k$ at $y$ is $\operatorname{Fitt}_{k-1}((N_y^\vee)^{\oplus m})$.

For a closed point $\alpha\in\widehat A$, put $y=\rho(\alpha)$. Since $\rho$ is smooth, $\widehat\OO_{\widehat A,\alpha}$ is flat over $R_y$. The completion isomorphism in~(i), together with Remark~\ref{rem:localFitting} and flat base change for degree-zero jump ideals, gives
\[
 \mathcal I_{Z_k,\alpha}\widehat\OO_{\widehat A,\alpha}
 =\operatorname{Fitt}_{k-1}\bigl((N_y^\vee)^{\oplus m}\bigr)\widehat\OO_{\widehat A,\alpha}
 =\mathcal I_{\rho^{-1}(D_k),\alpha}\widehat\OO_{\widehat A,\alpha}.
\]
Outside $\rho^{-1}(\operatorname{Supp}N)$ both ideals are the unit ideal. Faithful flatness of $\OO_{\widehat A,\alpha}\to\widehat\OO_{\widehat A,\alpha}$ yields equality of the two ideals in each local ring at a closed point. Since the ideal sheaves are coherent on the variety $\widehat A$, these equalities imply
\[
 Z_k=\rho^{-1}(D_k)
\]
scheme-theoretically. In particular, no gluing or algebraization of the formal module isomorphisms is required. Each nonempty finite subscheme $D_k$ has dimension zero, and $\rho$ is smooth of relative dimension $g-s$, proving the codimension assertion in~(ii).

Finally, a minimal free resolution of $N_y^{\oplus m}$ over $R_y$, shifted by $[-s]$, has degree-zero rank $m\dim_\CC\operatorname{soc}_{R_y}(N_y)$ by Remark~\ref{rem:localFitting}. After flat base change to $\widehat\OO_{\widehat A,\alpha}$ and then passage to the residue field, its differentials still vanish, because their entries lie in the maximal ideal of $R_y$. Derived base change thus gives~(iii). The two special cases follow by taking $K=A$ and $K=0$, respectively; in the latter case $N$ is the fiber of $\mathcal W\otimes L$ at the origin, a vector space of dimension $r$.
\end{proof}

\begin{remark}\label{rem:twistsum}
The Fitting ideals must be computed for the direct sum~\eqref{eq:Nsum}, not support point by support point. For instance, if two points of $\operatorname{Supp}M'$ lie over the same $y$ and each contributes a residue field, then $N_y=k^{\oplus2}$ and the ideal of $Z_k$ is $\mathfrak m_y^{\max(2m-k+1,0)}$, whereas intersecting the separate ideals would give $\mathfrak m_y^{m}$ for $k=1$. Similarly, the socle in~(iii) is taken over $R_y$: the $\CC[[x,y]]$-module $\CC[[x,y]]/(x^2,y)$ has a one-dimensional socle, but its restriction to $\CC[[y]]$ is $k^{\oplus2}$.
\end{remark}

\subsection{Semihomogeneous bundles}\label{subsec:semihom}

A vector bundle $E$ on $A$ is \emph{semihomogeneous} if for every $x\in A$ there is $P\in\widehat A$ with $t_x^*E\simeq E\otimes P$~\cite{Mukai78}. Put $\delta(E)=c_1(E)/\rk E\in\operatorname{NS}(A)_{\mathbb Q}$. Every bundle $\mathcal W\otimes L$ as in Theorem~\ref{thm:twistjump} is semihomogeneous, but not conversely (Example~\ref{ex:semihom}). We use the following results of Mukai~\cite{Mukai78}, valid over $\CC$:
\begin{enumerate}[label=(\alph*)]
\item every semihomogeneous bundle is isomorphic to $\bigoplus_iU_i\otimes E_i$, where each $E_i$ is a simple semihomogeneous bundle with $\delta(E_i)=\delta(E)$ and each $U_i$ is unipotent, that is, a successive extension of copies of $\OO_A$ (Propositions~6.18 and~6.2);
\item two simple semihomogeneous bundles with the same $\delta$ differ by a twist by an element of $\widehat A$ (Proposition~6.17(1));
\item a simple semihomogeneous bundle is isomorphic to $f_*L'$ for an isogeny $f:Y\to A$ and a line bundle $L'$ on $Y$ (Theorem~5.8);
\item for an isogeny $f$ and a vector bundle $F$ on $Y$, $f^*f_*F\simeq\bigoplus_{x\in\ker f}t_x^*F$ (proof of Lemma~6.21); indeed $Y\times_AY$ is the disjoint union of the graphs of the translations by $\ker f$.
\end{enumerate}

\begin{lemma}\label{lem:semihomreduce}
Let $E$ be a semihomogeneous bundle of rank $R$ on $A$. There are an isogeny $f:Y\to A$, a homogeneous bundle $\mathcal W_Y$ on $Y$ and a line bundle $L'$ on $Y$ such that, with $G=\mathcal W_Y\otimes L'$,
\[
 E\simeq f_*G,\qquad f^*c_1(\det E)=R\,c_1(L'),\qquad \Phi_A(E)\simeq\widehat f^*\Phi_Y(G),
\]
where $\widehat f:\widehat A\to\widehat Y$ is the dual isogeny. Consequently $Z_k(E)=\widehat f^{-1}(Z_k(G))$ as schemes for every $k$. Moreover, $E$ is nef if and only if $L'$ is nef, if and only if $\delta(E)$ is nef.
\end{lemma}

\begin{proof}
By~(a) and~(b), fixing a summand $E_0$, each $E_i$ is isomorphic to $E_0\otimes P_{\alpha_i}$ for some $P_{\alpha_i}\in\widehat A$, so $E\simeq\mathcal W\otimes E_0$ where $\mathcal W:=\bigoplus_i(U_i\otimes P_{\alpha_i})$ is $\widehat A$-filtered and hence homogeneous. By~(c), $E_0\simeq f_*L'$ for an isogeny $f:Y\to A$ and a line bundle $L'$ on $Y$. The projection formula gives
\[
 E\simeq\mathcal W\otimes f_*L'\simeq f_*(f^*\mathcal W\otimes L')=f_*G,
\]
where $\mathcal W_Y:=f^*\mathcal W$ is $\widehat Y$-filtered (hence homogeneous) and $G=\mathcal W_Y\otimes L'$. By~(d), $f^*E\simeq f^*f_*G\simeq\bigoplus_{x\in\ker f}t_x^*G$. Since $c_1(\mathcal W_Y)=0$ and $R=\rk E=(\deg f)\rk\mathcal W_Y=|\ker f|\rk\mathcal W_Y$, taking first Chern classes yields
\[
 f^*c_1(\det E)=|\ker f|\,c_1(G)=R\,c_1(L').
\]

For the normalized Poincar\'e bundles, $(f\times1_{\widehat A})^*\mathcal P_A\simeq(1_Y\times\widehat f)^*\mathcal P_Y$~\cite[proof of (3.4)]{Mukai81}. By the projection formula,
\[
 \pr_A^*E\otimes\mathcal P_A\simeq(f\times1)_*\bigl((1\times\widehat f)^*(\pr_Y^*G\otimes\mathcal P_Y)\bigr).
\]
As $f\times1$ is finite and $\widehat f$ is flat, pushing forward to $\widehat A$ and using flat base change give $\Phi_A(E)\simeq\widehat f^*\Phi_Y(G)$, with no shift since $\dim Y=\dim A$. Degree-zero jump ideals commute with flat base change, which gives $Z_k(E)=\widehat f^{-1}(Z_k(G))$.

If $E$ is nef, so is $\det E$, hence $\delta(E)$ and $c_1(L')=R^{-1}f^*c_1(\det E)$ are nef. If $L'$ is nef, then by~(d) $f^*E\simeq\bigoplus_{x\in\ker f}f^*\mathcal W\otimes t_x^*L'$ has a filtration by nef line bundles, so $f^*E$ and therefore $E$ are nef~\cite[Theorem~6.2.12]{Lazarsfeld}.
\end{proof}

Applying Theorem~\ref{thm:twistjump} to $G$ on $Y$ and comparing the connected kernels $K(L')^0$ and $K(\det E)^0$ under the isogeny $f$ now yields the general codimension theorem for nef semihomogeneous bundles.

\begin{theorem}\label{thm:semihomjump}
Let $E$ be a nef semihomogeneous bundle on $A$, put $K_E=K(\det E)^0$, $s=\dim K_E$, and let $\rho_E:\widehat A\to\widehat K_E$ be restriction. There are finite subschemes $D_k\subset\widehat K_E$ with $Z_k(E)=\rho_E^{-1}(D_k)$ for all $k\ge1$; every nonempty $Z_k(E)\ne\widehat A$ has codimension exactly $s$. If $E=\mathcal W\otimes L$ as in Theorem~\ref{thm:twistjump}, then $K_E=K(L)^0$ and $D_k$ is the Fitting scheme of Theorem~\ref{thm:twistjump}(ii).
\end{theorem}

\begin{proof}
Take $f$ and $G=\mathcal W_Y\otimes L'$ as in Lemma~\ref{lem:semihomreduce}; $L'$ is nef. Put $K_Y=K(L')^0$. Theorem~\ref{thm:twistjump} on $Y$ gives finite subschemes $D'_k\subset\widehat K_Y$ with $Z_k(G)=\rho_Y^{-1}(D'_k)$. Because $f^*c_1(\det E)=R\,c_1(L')$ and $\phi_M$ depends only on the numerical class of $M$, we have $\widehat f\circ\phi_{\det E}\circ f=R\phi_{L'}$. Hence $\phi_{\det E}\circ f$ maps the connected group $K_Y$ into the finite group $\ker\widehat f$, hence to zero, so $f(K_Y)\subset K_E$. For the associated translation-invariant Hermitian forms one has $H_{\det E}(df(v),df(w))=R\,H_{L'}(v,w)$; since $df$ is a complex linear isomorphism and $R>0$, it maps $T_0K_Y=\ker H_{L'}$ onto $T_0K_E=\ker H_{\det E}$. Thus $\dim K_Y=\dim K_E$, so $f(K_Y)=K_E$, and $\lambda:=f|_{K_Y}:K_Y\to K_E$ is an isogeny, its kernel lying in $\ker f$. Both $\rho_Y\circ\widehat f$ and $\widehat\lambda\circ\rho_E$ send $P\in\widehat A$ to $(f^*P)|_{K_Y}$, so
\[
 Z_k(E)=\widehat f^{-1}\rho_Y^{-1}(D'_k)=\rho_E^{-1}\bigl(\widehat\lambda^{-1}(D'_k)\bigr).
\]
Put $D_k:=\widehat\lambda^{-1}(D'_k)$, which is finite since $\widehat\lambda$ is an isogeny. Since $\rho_E$ is smooth and surjective with fibers of dimension $g-s$, every nonempty $Z_k(E)\ne\widehat A$ has codimension $s$ in $\widehat A$. The final assertion is Theorem~\ref{thm:twistjump}(ii).
\end{proof}

\subsection{Elliptic directions and abelian surfaces}\label{subsec:surfacejump}

We now specialize to the case $s=1$ and to abelian surfaces, which are the cases used in Section~\ref{sec:classification}. When $s=1$, so that $K$ is an elliptic curve, the local modules are explicit. On an elliptic curve $C$, let $F_r$ denote Atiyah's indecomposable unipotent bundle of rank $r$~\cite{Atiyah}. A homogeneous bundle $\mathcal U$ on $C$ is a direct sum $\bigoplus_jN_j\otimes F_{r_j}$ with $N_j\in\widehat C$. For $y\in\widehat C$, the \emph{type of $\mathcal U$ at $y$} is the multiset $\{r_j:N_j\otimes P_y\simeq\OO_C\}$; it is nonempty for only finitely many $y$.

\begin{corollary}\label{cor:ellipticjump}
In Theorem~\ref{thm:twistjump}, suppose $s=1$, so that $C:=K$ is an elliptic curve. For $y\in\widehat C$, repeat every part of the type of $(\mathcal W\otimes L)|_C$ at $y$ exactly $m$ times. If the resulting list has length $n$, then near $\tau_y$ the scheme $Z_k$ is the Cartier divisor $e_k\,\tau_y$, where $e_k$ is the sum of the $n-k+1$ smallest entries for $k\le n$, and $e_k=0$ otherwise. In particular, if $g=2$, $L$ is nef, $c=c_1(L)\not\equiv0$ and $c^2=0$, then $s=1$ and $c\equiv m[C]$.
\end{corollary}

\begin{proof}
On $C$, the transform $\Phi_C(F_r)\simeq M_{F_r}[-1]$ is a length-$r$ sheaf supported at $0\in\widehat C$ by its filtration. Derived base change to $0\in\widehat C$ identifies $M_{F_r}\otimes k(0)$ with $H^1(C,F_r)$ of dimension $1$, so $M_{F_r}$ is cyclic by Nakayama's lemma and hence isomorphic to $\CC[[y]]/(y^r)$ over the discrete valuation ring $R_0\simeq\CC[[y]]$; twisting by a character translates its support. For each cyclic summand $\CC[[y]]/(y^r)$, dualizing the resolution $0\to\CC[[y]]\xrightarrow{y^r}\CC[[y]]\to\CC[[y]]/(y^r)\to0$ gives $\Ext^1_{\CC[[y]]}(\CC[[y]]/(y^r),\CC[[y]])\simeq\CC[[y]]/(y^r)$. Hence $N_y^\vee\simeq N_y$, and $(N_y^\vee)^{\oplus m}\simeq\bigoplus_{i=1}^n\CC[[y]]/(y^{r_i})$ with $r_1\le\cdots\le r_n$. Its minimal presentation matrix is $\operatorname{diag}(y^{r_1},\dots,y^{r_n})$, whose minors of order $n-k+1$ generate $(y^{e_k})$ for $1\le k\le n$ (and the unit ideal for $k>n$). Since $\rho$ is smooth, pulling back $(y^{e_k})$ along $\rho$ defines $e_k\tau_y$, giving the first assertion by Theorem~\ref{thm:twistjump}(ii).

If $g=2$ and $c^2=0\ne c$, then $\chi(L)=c^2/2=0$, so $\phi_L$ cannot be an isogeny (as $\deg\phi_L=\chi(L)^2=0$) and is nonzero (as $c\not\equiv0$); hence $s=1$, $\bar A$ is an elliptic curve, and $m=\deg M$. Since the pullback under $f$ of a point class of $\bar A$ is the class of a fiber of $f$ (a translate of $C$), the descent~\eqref{eq:nefdescent} gives $c\equiv m[C]$.
\end{proof}

On an abelian surface, the dimension $s=\dim K_E$ can only be $2$, $1$, or $0$, determined by the numerical class and self-intersection of $\delta(E)=c_1(E)/\rk E$. Combining Theorem~\ref{thm:semihomjump} with Corollary~\ref{cor:ellipticjump} yields the following trichotomy.

\begin{corollary}\label{cor:semihomsurface}
Let $A$ be an abelian surface and $E$ a nef semihomogeneous bundle with $\delta=\delta(E)$. Then $Z_k(E)$ is finite if $\delta\equiv0$; an effective Cartier divisor supported on finitely many translates of an elliptic curve if $\delta\not\equiv0$ and $\delta^2=0$; and $\widehat A$ or empty if $\delta$ is ample. For $E=\mathcal W\otimes L$ the multiplicities in the second case are those of Corollary~\ref{cor:ellipticjump}.
\end{corollary}

\begin{proof}
The Hermitian form $H_\delta$ of the nef class $\delta$ is positive semidefinite, and $s=\dim\ker H_\delta$. If $\lambda_1,\lambda_2\ge0$ are its eigenvalues with respect to a fixed positive definite form, then $\delta^2=C\lambda_1\lambda_2$ for some $C>0$. Hence $H_\delta=0$ if and only if $\delta\equiv0$, and $H_\delta$ is nondegenerate if and only if $\delta^2>0$, that is, if and only if $\delta$ is ample. Thus $s=2,1,0$ in the three cases. For $s=2$, $\rho_E$ is the identity, so $Z_k(E)=D_k$ is finite by Theorem~\ref{thm:semihomjump}. For $s=1$, $\widehat K_E$ is an elliptic curve, so any finite subscheme $D_k\subset\widehat K_E$ is an effective Cartier divisor, and its preimage under the smooth morphism $\rho_E$ (whose fibers are translates of an elliptic curve) is again an effective Cartier divisor supported on translates of an elliptic curve. For $s=0$, $\widehat K_E$ is a point, so $D_k$ is either that point or empty, giving $Z_k(E)=\widehat A$ or $\varnothing$. The last sentence follows from Corollary~\ref{cor:ellipticjump}.
\end{proof}

\begin{example}\label{ex:semihom}
Let $A=C_1\times C_2$ be a product of elliptic curves, let $F$ be a stable bundle of rank two and degree one on $C_1$, and put $E=p_1^*F$. It is semihomogeneous, since $t_x^*F\simeq F\otimes P$ for suitable $P$ on $C_1$, but it is not of the form $\mathcal W\otimes L$: for such a bundle, $\det$ has even degree on the curves $C_1\times\{c\}$. Here $\delta(E)=\tfrac12p_1^*[\mathrm{pt}]$ is nef with $\delta^2=0$, and $K_E=\{0\}\times C_2$. Since $h^0(F\otimes P_1)=1$ and $h^1(F\otimes P_1)=0$ for all $P_1\in\widehat C_1$, the K\"unneth formula gives the local complex $[\OO\xrightarrow{t}\OO]$ along $\widehat C_1\times\{0\}$, so $Z_1(E)=\widehat C_1\times\{0\}$ is reduced and $Z_2(E)=\varnothing$, in accordance with Corollary~\ref{cor:semihomsurface}.
\end{example}

\section{The anticanonical trichotomy and its examples}\label{sec:classification}

We now return to the projective bundle $X$ and prove Theorem~\ref{thm:introC}. First we show that $-K_X$ is nef exactly when $\Sym^2\VV$ is a homogeneous bundle twisted by $\det\VV$, so that Theorem~\ref{thm:twistjump} applies. Section~\ref{subsec:homogscroll} then settles the homogeneous case $c\equiv0$ (Theorem~\ref{thm:homogeneous}), and Section~\ref{subsec:trichotomy} combines it with Theorem~\ref{thm:twistjump} to obtain the trichotomy (Corollary~\ref{cor:trichotomy}). Section~\ref{subsec:examples} gives examples of the second case and shows that the hypothesis on $-K_X$ cannot be dropped.

Throughout this section $\VV$ is nef, $c=c_1(\VV)$, and $a$ satisfies \eqref{eq:oddvanish} and \eqref{eq:oddachoice}, so that $\mathcal D_a$ exists by Theorem~\ref{thm:oddcomponent}. By~\eqref{eq:oddBNtransfer}, $\BN_X^k\cap\mathcal D_a=q_2^{-1}(Z_k)$, where $Z_k\subset\widehat A$ is the degree-zero jump scheme of the family $\{\Sym^2\VV\otimes P\}_{P\in\widehat A}$ for $h^0\ge k$. This family on $A$ is independent of $a$. Put
\[
 \mathcal G=\Sym^2\VV\otimes(\det\VV)^{-1},
\]
a rank-three bundle, so that $\Sym^2\VV\otimes P=\mathcal G\otimes\det\VV\otimes P$.

\subsection{The anticanonical class and homogeneous bundles}\label{subsec:homogscroll}

We first show that $-K_X$ is nef exactly when $\mathcal G$ is homogeneous, and then treat completely the basic case $c\equiv0$, in which $\VV$ itself is homogeneous.

\begin{lemma}\label{lem:antinef}
For $X=\PP_A(\VV)$ and $\mathcal G=\Sym^2\VV\otimes(\det\VV)^{-1}$, the following conditions are equivalent:
\begin{enumerate}[label=(\alph*)]
\item the anticanonical class $-K_X$ is nef;
\item the rank-three bundle $\mathcal G$ is numerically flat;
\item $\mathcal G$ is $\widehat A$-filtered (equivalently, homogeneous).
\end{enumerate}
\end{lemma}

\begin{proof}
Since $-K_X=2\bigl(\zeta-\tfrac12\pi^*c\bigr)$, condition~(a) means that the $\mathbb Q$-twisted bundle $\VV\bigl\langle-\tfrac12 c\bigr\rangle$ is nef~\cite[Section~6.2]{Lazarsfeld}. Because $\VV$ has rank two, $\VV^*\simeq\VV\otimes(\det\VV)^{-1}$, so the dual of $\VV\bigl\langle-\tfrac12 c\bigr\rangle$ is
\[
 \Bigl(\VV\bigl\langle-\tfrac12 c\bigr\rangle\Bigr)^*
 =\VV^*\bigl\langle\tfrac12 c\bigr\rangle
 \simeq\VV\bigl\langle-c+\tfrac12 c\bigr\rangle
 =\VV\bigl\langle-\tfrac12 c\bigr\rangle.
\]
Hence~(a) is equivalent to numerical flatness of $\VV\bigl\langle-\tfrac12 c\bigr\rangle$. Symmetric powers of nef $\mathbb Q$-twisted bundles are nef~\cite[Theorem~6.2.12]{Lazarsfeld}, so $\mathcal G=\Sym^2\bigl(\VV\bigl\langle-\tfrac12 c\bigr\rangle\bigr)$ is nef; moreover,
\[
 \mathcal G^*=\Sym^2(\VV^*)\otimes\det\VV\simeq\Sym^2\bigl(\VV\otimes(\det\VV)^{-1}\bigr)\otimes\det\VV\simeq\mathcal G,
\]
so $\mathcal G^*$ is also nef and $\mathcal G$ is numerically flat. Thus~(a) implies~(b). Conversely, under the relative Veronese embedding $v:X=\PP_A(\VV)\hookrightarrow\PP_A(\Sym^2\VV)\simeq\PP_A(\mathcal G)$, the tautological class $\zeta_{\mathcal G}$ pulls back to $2\zeta-\pi^*c=-K_X$; if $\mathcal G$ is nef, then $\zeta_{\mathcal G}$ is nef and its restriction $-K_X$ to $X$ is nef.

By Demailly--Peternell--Schneider~\cite[Theorem~1.18]{DPS}, a numerically flat bundle admits a filtration by subbundles whose graded pieces are Hermitian flat, that is, given by unitary representations of $\pi_1(A)\simeq\ZZ^4$. Unitary representations of an abelian group are direct sums of characters, and flat unitary line bundles on $A$ are precisely the elements of $\widehat A$. Refining the filtration gives~(c). Conversely, an $\widehat A$-filtered bundle is a successive extension of numerically flat line bundles and hence numerically flat~\cite[Theorem~6.2.12]{Lazarsfeld}, giving~(c)$\Rightarrow$(b).
\end{proof}

In particular, when $c\equiv0$, the nefness of $\VV$ implies that $\VV^*\simeq\VV\otimes(\det\VV)^{-1}$ is also nef, so $\VV$ itself is numerically flat and hence homogeneous. We now classify all resulting jump schemes in this basic case. For $0\ne\eta\in H^1(A,\OO_A)$, write $\mathcal U_\eta$ for the extension
\begin{equation}\label{eq:unipotent}
 0\to\OO_A\to\mathcal U_\eta\to\OO_A\to0.
\end{equation}
\begin{theorem}\label{thm:homogeneous}
Let $\VV$ be a rank-two $\widehat A$-filtered vector bundle on $(A,H)$ with $d=H^2$.
\begin{enumerate}[label=(\roman*)]
\item Up to twisting $\VV$ by a line bundle in $\widehat A$ (which leaves $X$ unchanged and translates the second Picard factor), either $\VV\simeq\OO_A\oplus Q$ for some $Q\in\widehat A$, or $\VV\simeq\mathcal U_\eta$ with $\eta\ne0$.
\item For every integer $a\ge1$ and rational $\tfrac{2a}{3}<t<2a$, the nonzero extensions~\eqref{eq:oddextension} form the smooth projective component $\mathcal D_a$ of Theorem~\ref{thm:oddcomponent}, which is a $\PP^{a^2d-1}$-bundle over $\widehat A^2$ of dimension $a^2d+3$. This interval is optimal by Corollary~\ref{cor:nefoptimal}, since the last graded piece of an $\widehat A$-filtered bundle is a numerically trivial quotient line bundle.
\item Scheme-theoretically, $\BN_X^k\cap\mathcal D_a=q_2^{-1}(Z_k)$, where $Z_k\subset\widehat A$ is the closed subscheme defined by the Fitting ideal $\operatorname{Fitt}_{k-1}$ of $\mathcal{E}xt^2_{\widehat A}(M_{\Sym^2\VV},\OO_{\widehat A})$. Denoting by $\mathfrak m_y$ the maximal ideal of $y\in\widehat A$ and by $[N]\in\widehat A$ the point corresponding to $N\in\Pic^0(A)$, we have $Z_k=\varnothing$ for all $k\ge4$, and up to the translation induced by a twist in~(i), for $1\le k\le3$:
\[
\renewcommand{\arraystretch}{1.2}\setlength{\arraycolsep}{4pt}
\begin{array}{@{}llll@{}}
\toprule
 \VV & \II_{Z_1} & \II_{Z_2} & \II_{Z_3}\\\midrule
 \OO_A^{\oplus2} & \mathfrak m_0^3 & \mathfrak m_0^2 & \mathfrak m_0\\
 \OO_A\oplus Q,\ \operatorname{ord}Q=2\, & \mathfrak m_0^2\cap\mathfrak m_{[Q]} & \mathfrak m_0 & (1)\\
 \OO_A\oplus Q,\ \operatorname{ord}Q\ge3\, & \mathfrak m_0\cap\mathfrak m_{[Q^{-1}]}\cap\mathfrak m_{[Q^{-2}]} & (1) & (1)\\
 \mathcal U_\eta & (y,x^3)\ \text{at }0 & (1) & (1)\\
\bottomrule
\end{array}
\]
Here $(1)$ is the unit ideal, $\operatorname{ord}Q\ge3$ includes infinite order, and in the last row $x,y$ are formal coordinates at the origin.

\end{enumerate}
In particular, every nonempty Brill--Noether scheme on $\mathcal D_a$ when $\VV$ is homogeneous is pulled back from a finite subscheme of $\widehat A$ of length at most $6$.
\end{theorem}

\begin{proof}
(i) Let $0\to N_1\to\VV\to N_2\to0$ with $N_i\in\widehat A$. If the extension splits, then $\VV\simeq N_2\otimes(\OO_A\oplus N_1N_2^{-1})$. Otherwise $H^1(A,N_1N_2^{-1})\ne0$, which forces $N_1\simeq N_2$ since nontrivial elements of $\widehat A$ have vanishing cohomology, and then $\VV\simeq N_1\otimes\mathcal U_\eta$ with $\eta\ne0$. Replacing $\VV$ by $\VV\otimes N$ replaces $\zeta$ by $\zeta+\pi^*c_1(N)$ and turns~\eqref{eq:oddextension} into the same family with $(P_1,P_2)$ replaced by $(P_1N^{-1},P_2N^2)$.

(ii) The bundles $W(aH)\otimes P$ and $\Sym^3\VV(-aH)\otimes P$ are filtered by line bundles numerically equivalent to $aH$ and $-aH$, respectively. By Kodaira vanishing and Riemann--Roch on $A$, \eqref{eq:oddvanish} holds for every $a\ge1$, with $h=\chi(W(aH))=2\chi(\OO_A(aH))=a^2d$. Since $c_1(\VV)\equiv0$ and $c_2(\VV)=0$, we have $A_t=t^2d$, $h_t=2td$, and $c_t=0$. Thus $aH-3c_1(\VV)$ is ample, \eqref{eq:oddachoice} holds, and~\eqref{eq:nefcommon} reads $\tfrac{2a}{3}<t<2a$, so Theorem~\ref{thm:oddcomponent} applies.

(iii) By~\eqref{eq:oddBNtransfer} and Proposition~\ref{prop:FMfitting} (with $g=2$), the ideal defining $Z_k$ at a point $y\in\widehat A$ is $\operatorname{Fitt}_{k-1}\bigl(\Ext^2_{R_y}(M_{\Sym^2\VV,y},R_y)\bigr)$, and $Z_k=\varnothing$ for $k>\rk\Sym^2\VV=3$.

If $\VV\simeq\OO_A\oplus Q$, then $\Sym^2\VV\simeq\OO_A\oplus Q\oplus Q^2$, so $M_{\Sym^2\VV}\simeq k(0)\oplus k([Q^{-1}])\oplus k([Q^{-2}])$ by Proposition~\ref{prop:FMfitting}. Over the two-dimensional regular local ring $R_y$ with maximal ideal $\mathfrak m_y=(x,y)$, the residue field $k(y)$ is Gorenstein ($\Ext^2_{R_y}(k(y),R_y)\simeq k(y)$) and presented by the $1\times2$ row matrix $(x,y)$. Hence if $m\in\{1,2,3\}$ summands coincide at $y$, then $\Ext^2_{R_y}(k(y)^{\oplus m},R_y)\simeq k(y)^{\oplus m}$ has a block-diagonal presentation matrix of size $m\times 2m$, whose minors of order $m-k+1$ generate $\mathfrak m_y^{m-k+1}$ for $1\le k\le m$ and $(1)$ for $k>m$. Distinguishing $Q=\OO_A$ ($m=3$ at $0$), $\operatorname{ord}Q=2$ ($m=2$ at $0$ and $m=1$ at $[Q]$), and $\operatorname{ord}Q\ge3$ ($m=1$ at $0,[Q^{-1}],[Q^{-2}]$) yields the first three rows of the table.

Now let $\VV=\mathcal U_\eta$ with $\eta\ne0$. Multiplication by the distinguished section of $\mathcal U_\eta$ yields the exact sequence
\[
 0\longrightarrow\mathcal U_\eta\longrightarrow\Sym^2\mathcal U_\eta\longrightarrow\OO_A\longrightarrow 0.
\]
In a local frame $\{e_0,e_1\}$ where the transition of $\mathcal U_\eta$ is $e_1\mapsto e_1+\eta_{ij}e_0$, the section $e_1^2$ transforms to $e_1^2+2\eta_{ij}e_0e_1+\eta_{ij}^2e_0^2$, so the connecting map $H^0(\OO_A)\to H^1(\mathcal U_\eta)$ composed with $H^1(\mathcal U_\eta)\to H^1(\OO_A)$ sends $1$ to $2\eta\ne0$. Since $h^0(\mathcal U_\eta)=1$ by~\eqref{eq:unipotent}, this forces $h^0(\Sym^2\mathcal U_\eta)=1$. As $\det\mathcal U_\eta=\OO_A$, both $\mathcal U_\eta$ and $\Sym^2\mathcal U_\eta$ are self-dual, so Serre duality gives $h^2(\Sym^2\mathcal U_\eta)=h^0(\Sym^2\mathcal U_\eta)=1$.

By Proposition~\ref{prop:FMfitting}, $\Phi_A(\Sym^2\mathcal U_\eta)\simeq\mathcal T_\eta[-2]$ for a sheaf $\mathcal T_\eta$ of length $3$ supported at $0\in\widehat A$. Derived base change gives $\dim_\CC(\mathcal T_\eta\otimes k(0))=h^2(\Sym^2\mathcal U_\eta)=1$, so $\mathcal T_\eta\simeq\widehat R_0/I$ is cyclic of length $3$ by Nakayama's lemma. Up to isomorphism, the only length-$3$ quotients of $\widehat R_0\simeq\CC[[x,y]]$ are $\CC[[x,y]]/\mathfrak m_0^2$ (socle dimension $2$) and $\CC[[x,y]]/(y,x^3)$ (socle dimension $1$). By Proposition~\ref{prop:FMfitting}, $\dim_\CC\operatorname{soc}(\mathcal T_\eta)=h^0(\Sym^2\mathcal U_\eta)=1$, so after a formal change of coordinates we may take $I=(y,x^3)$. Since $\CC[[x,y]]/(y,x^3)$ is a complete intersection (hence Gorenstein), its local $\Ext^2$ dual is isomorphic to itself, giving $\mathcal I_{Z_1,0}=(y,x^3)$ and $Z_2=Z_3=\varnothing$. Finally, the maximum length among all four cases is $\dim_\CC\CC[[x,y]]/\mathfrak m_0^3=6$.

The shifted Koszul resolution of this complete intersection is
\begingroup
\setlength{\abovedisplayskip}{4pt plus 1pt minus 1pt}
\setlength{\belowdisplayskip}{4pt plus 1pt minus 1pt}
\setlength{\abovedisplayshortskip}{0pt}
\setlength{\belowdisplayshortskip}{4pt plus 1pt minus 1pt}
\[
 \OO\xrightarrow{(y,x^3)^{\mathsf T}}\OO^2
 \xrightarrow{(-x^3,y)}\OO
\]
\endgroup
in degrees $0,1,2$, so its degree-zero jump ideal for $k=1$ is $(y,x^3)$ and $Z_k=\varnothing$ for $k\ge2$. Derived restriction to the origin gives cohomology dimensions $(1,2,1)$, while outside the origin the complex is acyclic. Smooth pullback by $q_2$ therefore gives a local complete intersection of codimension two, of transverse length three and support dimension $a^2d+1$, along which the cohomology dimensions on $X$ are $(1,2,1,0)$; all cohomology vanishes off it. Since completion is faithfully flat, the local complete-intersection and multiplicity assertions descend to the algebraic local rings.
\end{proof}

\begin{example}
Take $H^2=2$ and $a=1$. The polarization $\zeta+B$ yields a smooth projective five-dimensional component on $\PP_A(\mathcal U_\eta)$, with
\[
 (c_1,c_2,c_3)=(\zeta+B,2\zeta B,0).
\]
Its first Brill--Noether scheme is supported on a smooth threefold, a $\PP^1$-bundle over $\widehat A$, and has transverse algebra $\CC[x]/(x^3)$. Its Zariski tangent space along that support has dimension four.
\end{example}

For the trivial bundle $\VV=\OO_A^{\oplus2}$ and $a=1$, the first row gives $\II_{\BN_X^k\cap\mathcal D_1}=\II_Y^{4-k}$ for $1\le k\le3$, where $Y=q_2^{-1}(0)$. Thus $(h^0,h^1,h^2,h^3)=(3,6,3,0)$ on $Y$ and all cohomology vanishes elsewhere.
\begin{example}
Let $H$ be a principal polarization on an abelian surface, so $H^2=2$. With polarization $\zeta+B$, the bundles \eqref{eq:oddproductextension} form a smooth projective five-dimensional component, a $\PP^1$-bundle over $\widehat A^2$. Its $\BN^3$ is a smooth threefold, a $\PP^1$-bundle over $\widehat A$. The schemes $\BN^2$ and $\BN^1$ are its square-ideal and cube-ideal thickenings. Fixing $\det E=\OO_X(1)\otimes\pi^*H$ gives a three-dimensional component whose $\BN^3$ support is a single $\PP^1$; the lower levels retain the stated transverse power ideals.
\end{example}

\subsection{The anticanonical trichotomy}\label{subsec:trichotomy}

With the homogeneous case in hand, Theorem~\ref{thm:twistjump} gives the complete classification when $-K_X$ is nef.

\begin{corollary}\label{cor:trichotomy}
Let $\VV$ be a nef rank-two vector bundle on $A$ such that $-K_X$ is nef, put $c=c_1(\VV)$, and let $Z_k\subset\widehat A$ be the jump schemes of $\{\Sym^2\VV\otimes P\}_{P\in\widehat A}$.
\begin{enumerate}[label=(\roman*)]
\item If $c\equiv0$, then $\VV$ is homogeneous, $\mathcal D_a$ exists for every $a\ge1$ and $\tfrac{2a}3<t<2a$, and the $Z_k$ are given by Theorem~\ref{thm:homogeneous}.
\item If $c\not\equiv0$ and $c^2=0$, let $C$, $m$, $\rho$ be as in Corollary~\ref{cor:ellipticjump} for $L=\det\VV$, and put $\delta=(\det\VV|_C)^{-1}\in\widehat C$. Then $\VV|_C$ is homogeneous and, up to a twist in $\widehat C$, one of the bundles in the left column below; correspondingly,
\[
\renewcommand{\arraystretch}{1.25}\setlength{\arraycolsep}{4pt}
\begin{array}{@{}ll@{}}
\toprule
 \VV|_C\ \text{up to twist} & Z_k\\\midrule
 \OO_C^{\oplus2} &(3m-k+1)\,\tau_{\delta}\\
 \OO_C\oplus R,\ \operatorname{ord}R=2 &(2m-k+1)\,\tau_{\delta R}+(m-k+1)\,\tau_{\delta}\\
 \OO_C\oplus R,\ \operatorname{ord}R\ge3 &(m-k+1)\,(\tau_{\delta}+\tau_{\delta R}+\tau_{\delta R^{-1}})\\
 F_2 &3(m-k+1)\,\tau_{\delta}\\
\bottomrule
\end{array}
\]
where $\tau_y=\rho^{-1}(y)$, points of $\widehat C$ are written as line bundles, $F_2$ is the nonsplit self-extension of $\OO_C$, $\operatorname{ord}R\ge3$ includes infinite order, and negative coefficients are read as zero. Every nonempty $Z_k$ has codimension one, whereas its two-term expected codimension is $k^2$.
\item If $c$ is ample, then $Z_k=\widehat A$ for $1\le k\le\tfrac32c^2$ and $Z_k=\varnothing$ otherwise, by the last assertion of Theorem~\ref{thm:twistjump} with $r=3$.
\end{enumerate}
\end{corollary}

\begin{proof}
Lemma~\ref{lem:antinef} gives $\Sym^2\VV=\mathcal G\otimes\det\VV$ with $\mathcal G$ homogeneous of rank three, so Theorem~\ref{thm:twistjump} applies with $\mathcal W=\mathcal G$ and $L=\det\VV$; this gives~(iii). If $c\equiv0$, then $\VV^*\simeq\VV\otimes(\det\VV)^{-1}$ is also nef, so $\VV$ is numerically flat and Theorem~\ref{thm:homogeneous} gives~(i).

In~(ii), Corollary~\ref{cor:ellipticjump} reduces the calculation to the Atiyah type of $(\mathcal G\otimes\det\VV)|_C$ along $C$. Since $c\cdot C=m[C]^2=0$ and $\VV$ is nef, $\VV|_C$ is a nef rank-two bundle of degree zero on the elliptic curve $C$, hence numerically flat and homogeneous. By Atiyah's classification~\cite{Atiyah}, there exists $B\in\widehat C$ such that either $\VV|_C\simeq B\otimes(\OO_C\oplus R)$ with $R\in\widehat C$, or $\VV|_C\simeq B\otimes F_2$. Since $\mathcal G=\Sym^2\VV\otimes(\det\VV)^{-1}$ is invariant under twisting $\VV$ by a line bundle, we have $\mathcal G|_C\simeq R\oplus\OO_C\oplus R^{-1}$ in the first case and $\mathcal G|_C\simeq\Sym^2F_2\simeq F_3$ in characteristic zero~\cite{Atiyah} in the second. Distinguishing $R=\OO_C$, $\operatorname{ord}R=2$, and $\operatorname{ord}R\ge3$ yields the types $(1,1,1)$, $(1,1)+(1)$, $(1)+(1)+(1)$, and $(3)$. A summand $N\otimes F_r$ of $\mathcal G|_C$ becomes $N\otimes\delta^{-1}\otimes F_r$ in $(\mathcal G\otimes\det\VV)|_C=\mathcal G|_C\otimes\delta^{-1}$ and therefore contributes at $N^{-1}\otimes\delta\in\widehat C$; the multiplicity formula of Corollary~\ref{cor:ellipticjump} then gives the table. Finally, by Theorem~\ref{thm:twistjump}(i) with $s=1$, all derived fibers of $\Phi_A(\Sym^2\VV)$ have cohomology concentrated in degrees $[0,1]$ (indeed, the DVR resolutions of Corollary~\ref{cor:ellipticjump} pull back locally to a two-term free complex $K^0\to K^1$ of equal rank, as $\chi=0$). The determinantal expected codimension of $\{h^0\ge k\}$ is therefore $k(k-\chi)=k^2$, whereas every nonempty $Z_k$ is an effective Cartier divisor and hence has codimension $1$.
\end{proof}

\begin{remark}
Since $-K_X$ is unchanged by twisting $\VV$, the three cases of Corollary~\ref{cor:trichotomy} depend on the chosen nef normalization of $\VV$, that is, on the chosen extension family~\eqref{eq:oddextension}. Twisting $\VV$ by a sufficiently ample line bundle always places the family in case~(iii).
\end{remark}

\begin{remark}[Projective bundles over abelian varieties of any dimension]\label{rem:anydim}
Let $A$ be an abelian variety of dimension $g$, let $\VV$ be a nef rank-two bundle on $A$ with $-K_X$ nef for $X=\PP_A(\VV)$, and let $Z_k\subset\widehat A$ be the jump schemes of $\{\Sym^2\VV\otimes P\}$. The proof of Lemma~\ref{lem:antinef} does not use $\dim A$, so $\mathcal G$ is homogeneous and Theorem~\ref{thm:twistjump} applies with $L=\det\VV$. Put $K=K(\det\VV)^0$ and $s=\dim K$. Since $c_1(\VV)|_K\equiv0$ and $\VV|_K$ is nef, $\VV|_K$ is numerically flat, hence homogeneous: up to a twist it is $\OO_K\oplus Q$ or the extension $\mathcal U_\eta$ of $\OO_K$ by itself with $0\ne\eta\in H^1(K,\OO_K)$. For $\mathcal U_\eta$, in additive coordinates of the formal group of $\widehat K$ with $\eta$ along the first axis, $M_{\mathcal U_\eta}=R/(x_1^2,x_2,\dots,x_s)$. The convolution square has underlying space $\CC[a,b]/(a^2,b^2)$ with $x_1$ acting by $w=a+b$; since $w^2=2ab$ and $w^3=0$, it splits as $\langle1,w,ab\rangle\oplus\langle a-b\rangle\simeq R/J_\eta\oplus k$ with $J_\eta=(x_1^3,x_2,\dots,x_s)$. As $\mathcal U_\eta^{\otimes2}\simeq\Sym^2\mathcal U_\eta\oplus\OO_K$, Krull--Schmidt gives $M_{\Sym^2\mathcal U_\eta}\simeq R/J_\eta$. For $s\ge1$, each cyclic summand $k$ or $R/J_\eta$ is an Artinian complete intersection and hence isomorphic to its local $\Ext^s$-dual, and $\operatorname{Fitt}_j(k^{\oplus n})=\mathfrak m^{n-j}$, $\operatorname{Fitt}_j((R/J_\eta)^{\oplus n})=J_\eta^{\,n-j}$. Setting $y_0=\delta$, $y_1=\delta Q$, and $y_2=\delta Q^{-1}$ where $\delta=(\det\VV|_K)^{-1}\in\widehat K$, Theorem~\ref{thm:twistjump} yields, with negative exponents read as the unit ideal,
\[
\renewcommand{\arraystretch}{1.2}\setlength{\arraycolsep}{4pt}
\begin{array}{@{}ll@{}}
\toprule
 \VV|_K\ \text{up to twist} & \text{ideal of }Z_k\text{ on }\widehat K\\\midrule
 \OO_K^{\oplus2} & \mathfrak m_{y_0}^{\,3m-k+1}\\
 \OO_K\oplus Q,\ \operatorname{ord}Q=2 & \mathfrak m_{y_1}^{\,2m-k+1}\cap\mathfrak m_{y_0}^{\,m-k+1}\\
 \OO_K\oplus Q,\ \operatorname{ord}Q\ge3 & \mathfrak m_{y_0}^{\,m-k+1}\cap\mathfrak m_{y_1}^{\,m-k+1}\cap\mathfrak m_{y_2}^{\,m-k+1}\\
 \mathcal U_\eta & J_\eta^{\,m-k+1}\\
\bottomrule
\end{array}
\]
pulled back along $\rho:\widehat A\to\widehat K$ (for $\mathcal U_\eta$, $J_\eta$ is taken at $y_0=\delta$). When $s=0$, $\det\VV$ is ample and $Z_k$ is either $\widehat A$ or empty as in Theorem~\ref{thm:twistjump}. Every nonempty proper $Z_k$ has codimension exactly $s$. For $g=2$, $s=2$ and $m=1$ this is the table of Theorem~\ref{thm:homogeneous}, and for $s=1$ it is Corollary~\ref{cor:trichotomy}(ii).
\end{remark}

\subsection{Examples}\label{subsec:examples}

We give two examples of case~(ii) of Corollary~\ref{cor:trichotomy}, at the two ends of its table, and an example showing that the hypothesis on $-K_X$ cannot be dropped. The first is an explicit instance with $m=1$ and $\VV|_C$ trivial, together with an independent direct proof.

Let $A=C_1\times C_2$ be a product of elliptic curves, with projections $p_i$. Choose degree-one line bundles $H_i$ on $C_i$ and put $H=p_1^*H_1\otimes p_2^*H_2$. On $C_1$, take the nonzero extension
\begin{equation}\label{eq:ellipticbundle}
 0\longrightarrow\OO_{C_1}\longrightarrow G\longrightarrow H_1\longrightarrow0.
\end{equation}
It is unique up to scaling the extension class because $h^1(C_1,H_1^{-1})=1$. Put $\VV=p_1^*G$ and $X=\PP_A(\VV)$. This $\VV$ is non-split and nonhomogeneous: restriction to $C_1$ detects nonsplitting, and its determinant $p_1^*H_1$ is not invariant under all translations.

\begin{theorem}\label{thm:ellipticBN}
For $X=\PP_{C_1\times C_2}(p_1^*G)$ as constructed above, $a=4$, and any rational $t$ satisfying $2<t<\tfrac{5+\sqrt{33}}2$ (for instance, the polarization $\zeta+3B$), the nonzero extensions~\eqref{eq:oddextension} form a smooth projective connected component $\mathcal D_4$ of dimension $23$, isomorphic to a $\PP^{19}$-bundle over $\widehat A^2$. Let $Y=q_2^{-1}(\widehat C_1\times\{0\})\subset\mathcal D_4$, which is a smooth prime divisor of dimension $22$. Then $\BN_X^k\cap\mathcal D_4=\varnothing$ for all $k\ge4$, and for $k=1,2,3$ the ideal sheaves satisfy
\begin{equation}\label{eq:ellipticBN}
 \II_{\BN_X^1\cap\mathcal D_4}=\II_Y^3,\qquad
 \II_{\BN_X^2\cap\mathcal D_4}=\II_Y^2,\qquad
 \II_{\BN_X^3\cap\mathcal D_4}=\II_Y.
\end{equation}
Moreover, any bundle $E\in\mathcal D_4$ has cohomology dimensions $(h^0,h^1,h^2,h^3)=(3,3,0,0)$ along $Y$ and $(0,0,0,0)$ on $\mathcal D_4\setminus Y$.
\end{theorem}

\begin{proof}
The bundle $G$ is stable on $C_1$: a line subbundle of degree at least one would map nontrivially into $H_1$, so it would have degree exactly one and map isomorphically, splitting~\eqref{eq:ellipticbundle}. As an extension of the nef line bundles $\OO_{C_1}$ and $H_1$, $G$ is nef, and so is $\VV=p_1^*G$; moreover, $G\langle-\tfrac12H_1\rangle$ is semistable of degree zero on $C_1$, hence numerically flat, so $\VV\langle-\tfrac12c_1(\VV)\rangle$ is numerically flat and $-K_X$ is nef by Lemma~\ref{lem:antinef}. Symmetric powers of $G$ are semistable in characteristic zero~\cite{HL}, so $\Sym^2G$ has rank three, degree three and positive slope, whence $H^1(C_1,\Sym^2G\otimes P)=0$ and $h^0(C_1,\Sym^2G\otimes P)=3$ for every $P\in\widehat C_1$. Likewise, $G^*\otimes H_1^{a-1}$ has slope $a-3/2>0$ and $\Sym^3G\otimes H_1^{-a}$ has slope $3/2-a<0$ for $a\ge2$. Combined with the positive or negative degree of $H_2^{\pm a}$ on $C_2$, the K\"unneth formula gives the vanishings in~\eqref{eq:oddvanish} uniformly on $\widehat A$ for $a\ge2$, with extension dimension $h=(2a-3)a$.

Here $H^2=2$, $c_1(\VV)\cdot H=1$, and $c_1(\VV)^2=c_2(\VV)=0$. For $a=4$ one has $h=20$, and $aH-3c_1(\VV)=p_1^*H_1+4p_2^*H_2$ is ample. In the notation of Theorem~\ref{thm:nefexact}, $h_t=1+4t$, $c_t=2t$ and $A_t=2t^2+2t$, so~\eqref{eq:nefcommon} reads $6t^2-10t-4>0$ and $t^2-5t-2<0$, that is, $2<t<\tfrac{5+\sqrt{33}}2$. Theorem~\ref{thm:nefexact} and the argument of Theorem~\ref{thm:oddcomponent} give the asserted component.

To compute its Brill--Noether scheme, first take relative cohomology of $\Sym^2G\otimes P$ over $\widehat C_1$. It is a vector bundle $\mathcal B$ of rank three in degree zero. On $\widehat C_2$, the Poincar\'e transform of $\OO_{C_2}$ is $k(0)[-1]$ \cite{Mukai81}. Thus, on $\widehat A=\widehat C_1\times\widehat C_2$, the relative cohomology of $\Sym^2\VV$ is the external product of $\mathcal B$ with $k(0)[-1]$. Locally along $\widehat C_1\times\{0\}$ it has the two-term representative
\[
 \OO^{\oplus3}\xrightarrow{yI_3}\OO^{\oplus3}
\]
in degrees zero and one, where $y$ cuts out that divisor. The minors of sizes $3,2,1$ generate $(y^3),(y^2),(y)$ respectively. Smooth pullback by $q_2$ gives \eqref{eq:ellipticBN}. The same complex, or the K\"unneth formula, gives the cohomology dimensions. The support dimension is $23-1=22$.
\end{proof}

Thus the three nonempty jump schemes are effective Cartier divisors with multiplicities $3,2,1$ on the same support. Here higher cohomology in degrees at least two vanishes on the entire component, so the two-term Brill--Noether expected codimension is applicable and equals $k^2$, since $\chi(E)=0$. In particular $\BN^1$ has the expected codimension one but is generically nonreduced, while $\BN^3$ has actual codimension one instead of the expected nine.

While Theorem~\ref{thm:ellipticBN} realizes the split case in the first row of Corollary~\ref{cor:trichotomy}(ii), the next example realizes the indecomposable fourth row.

\begin{example}\label{ex:twocubic}
Let $A=C\times C'$ with projections $p_1,p_2$, let $M$ have degree one on $C'$, and take $\VV=p_1^*F_2\otimes p_2^*M$. The pullback $p_1^*F_2$ is numerically flat, so $\VV$ is nef and its normalized tautological class is nef. Thus $-K_X$ is nef. Since $\det\VV=p_2^*M^2$, we have $m=2$ and $\VV|_C\simeq F_2$. On the two factors, $\Sym^2\VV=F_3\boxtimes M^2$. The transform of $F_3$ is the length-three module $\CC[[y]]/(y^3)$ at the origin in degree one, while $M^2$ transforms to a rank-two bundle in degree zero (as $\deg(M^2)=2>0$, so $H^1(C',M^2\otimes P_2)=0$ and $h^0(C',M^2\otimes P_2)=2$ for all $P_2\in\widehat{C'}$). By the K\"unneth formula the local cohomology complex along $\tau=\{0\}\times\widehat{C'}$ is
\[
 \OO^{\oplus2}\xrightarrow{\operatorname{diag}(y^3,y^3)}\OO^{\oplus2},
\]
where $y=0$ cuts out $\tau$. Its determinant gives $Z_1=6\tau$, its entries give $Z_2=3\tau$, and $Z_k=\varnothing$ for $k\ge3$. In particular the second locus detects both copies of the cubic block.
\end{example}

Finally, we show that the nefness hypothesis on $-K_X$ in Theorem~\ref{thm:introC} cannot be omitted.

\begin{example}\label{rem:mixed}
Let $f:A\to C'$ be an elliptic quotient and take $\VV=\OO_A\oplus f^*M$ with $m=\deg M>0$. Then $\VV$ is nef, whereas the quotient $\VV\twoheadrightarrow\OO_A$ defines a section $\sigma:A\hookrightarrow X$ with $\sigma^*\OO_X(1)\simeq\OO_A$; along any curve $C_0\subset\sigma(A)$ dominating $C'$, we have $\zeta\cdot C_0=0$ and $\pi^*c_1(\VV)\cdot C_0=m\deg(C_0\to C')>0$, so $(-K_X)\cdot C_0<0$ and $-K_X$ is not nef. The symmetric square is $\OO_A\oplus f^*M\oplus f^*M^2$. Near the origin choose parameters $x,y$ with $y=0$ defining $\widehat{C'}$. The first summand gives the Koszul column $(x,y)^{\mathsf T}$, and the positive-degree summands give $3m$ diagonal entries $y$ by Theorem~\ref{thm:twistjump} and its proof. For $1\le k\le3m+1$ the minors of order $3m+2-k$ of this block matrix generate
\[
 \mathcal I_{Z_k,0}=y^{3m+1-k}(x,y).
\]
Indeed a nonzero minor either uses only diagonal entries or uses one entry of the Koszul column. Away from the origin only the diagonal block contributes. Thus for $k\le3m$ the divisorial part is $(3m+1-k)\widehat{C'}$ with an embedded point at the origin. The last nonempty scheme $Z_{3m+1}$ is the reduced origin. More generally, for a split nef bundle the three line summands of its symmetric square give a block diagonal complex: each minor ideal is the sum of products of the individual minor ideals with the prescribed total order.
\end{example}

\medskip
\noindent\textbf{Acknowledgments.}
The author is deeply grateful to his supervisor, Professor Rong Du, for his
insightful advice, valuable discussions, and continuous support during
the preparation of this paper.
The author is supported by the National Natural Science Foundation of China
(Grant No.~12471040) and the Science and Technology Commission of Shanghai
Municipality (Grant No.~22DZ2229014).

\begingroup
\small
\setlength{\parskip}{0pt}

\endgroup

\begin{thebibliography}{99}
\bibitem{AP}
N. Alvarado and G. Pareschi, \emph{Semihomogenous vector bundles, $\mathbb Q$-twisted sheaves, duality, and linear systems on abelian varieties}, Int. J. Math. 37 (2026), no. 13, article 2650070.
\href{https://arxiv.org/abs/2407.20646}{arXiv:2407.20646}.
\bibitem{Atiyah}
M. F. Atiyah, \emph{Vector bundles over an elliptic curve}, Proc. London Math. Soc. (3) 7 (1957), 414--452.
\href{https://doi.org/10.1112/plms/s3-7.1.414}{doi:10.1112/plms/s3-7.1.414}.

\bibitem{BL}
A. Bayer and C. Li, \emph{Brill--Noether theory for curves on generic abelian surfaces}, Pure Appl. Math. Q. 13 (2017), no. 1, 49--76.
\href{https://arxiv.org/abs/1704.03546}{arXiv:1704.03546}.

\bibitem{BiL}
Ch. Birkenhake and H. Lange, \emph{Complex Abelian Varieties}, 2nd ed., Grundlehren der mathematischen Wissenschaften 302, Springer, 2004.

\bibitem{BW}
N. Budur and B. Wang, \emph{Cohomology jump loci of differential graded Lie algebras}, Compos. Math. 151 (2015), 1499--1528.
\href{https://arxiv.org/abs/1309.2264}{arXiv:1309.2264}.

\bibitem{Chien}
J. Chien, \emph{Codimension of Brill--Noether loci on abelian surfaces and the first Segre invariant}, preprint (2026).
\href{https://arxiv.org/abs/2610.24251}{arXiv:2610.24251}.

\bibitem{CNY}
I. Coskun, H. Nuer, and K. Yoshioka, \emph{Weak Brill--Noether on abelian surfaces}, Selecta Math. (N.S.) 31 (2025), article 48.
\href{https://arxiv.org/abs/2408.06095}{arXiv:2408.06095};
\href{https://doi.org/10.1007/s00029-025-01044-7}{doi:10.1007/s00029-025-01044-7}.

\bibitem{CMR}
L. Costa and R. M. Mir\'o-Roig, \emph{Brill--Noether theory for moduli spaces of sheaves on algebraic varieties}, Forum Math. 22 (2010), 411--432.
\href{https://arxiv.org/abs/0807.3232}{arXiv:0807.3232}.

\bibitem{CMT}
L. Costa and I. Mac\'ias Tarr\'io, \emph{Moduli spaces of stable bundles on ruled 3-folds and Brill--Noether problems}, Ricerche Mat. 75 (2026), 2179--2202. \href{https://doi.org/10.1007/s11587-026-01123-7}{doi:10.1007/s11587-026-01123-7}.

\bibitem{DPS}
J.-P. Demailly, T. Peternell, and M. Schneider, \emph{Compact complex manifolds with numerically effective tangent bundles}, J. Algebraic Geom. 3 (1994), 295--345.

\bibitem{Gulbrandsen}
M. G. Gulbrandsen, \emph{Fourier--Mukai transforms of line bundles on derived equivalent abelian varieties}, Matematiche (Catania) 63 (2008), 123--137.

\bibitem{HL}
D. Huybrechts and M. Lehn, \emph{The Geometry of Moduli Spaces of Sheaves}, 2nd ed., Cambridge University Press, 2010.

\bibitem{Lange}
H. Lange, \emph{Universal families of extensions}, J. Algebra 83 (1983), 101--112.

\bibitem{Lazarsfeld}
R. Lazarsfeld, \emph{Positivity in Algebraic Geometry II}, Ergebnisse der Mathematik und ihrer Grenzgebiete 49, Springer, 2004.

\bibitem{Miyanishi}
M. Miyanishi, \emph{Some remarks on algebraic homogeneous vector bundles}, in: Number Theory, Algebraic Geometry and Commutative Algebra, in honor of Y.~Akizuki, Kinokuniya, Tokyo, 1973, 71--93.

\bibitem{Mukai78}
S. Mukai, \emph{Semi-homogeneous vector bundles on an Abelian variety}, J. Math. Kyoto Univ. 18 (1978), 239--272.

\bibitem{Mukai81}
S. Mukai, \emph{Duality between $D(X)$ and $D(\widehat X)$ with its application to Picard sheaves}, Nagoya Math. J. 81 (1981), 153--175.
\href{https://doi.org/10.1017/S002776300001922X}{doi:10.1017/S002776300001922X}.

\bibitem{Mumford}
D. Mumford, \emph{Abelian Varieties}, Tata Institute of Fundamental Research Studies in Mathematics 5, Oxford University Press, 1970.

\bibitem{Polishchuk}
A. Polishchuk, \emph{Abelian Varieties, Theta Functions and the Fourier Transform}, Cambridge Tracts in Math. 153, Cambridge University Press, 2003.

\end{thebibliography}
\end{document}